\documentclass[12pt]{amsart}
\usepackage{amsmath,amsthm,amsfonts,amssymb}
\usepackage[all]{xy}
\input{grcawol.sty}
\input grcalcx.sty
\usepackage{color}

\DeclareMathOperator{\Ext}{Ext} \DeclareMathOperator{\Gr}{Gr}
\DeclareMathOperator{\Hom}{Hom} \DeclareMathOperator{\im}{im}
\DeclareMathOperator{\id}{id} \DeclareMathOperator{\ad}{ad}
\DeclareMathOperator{\Ht}{ht} \DeclareMathOperator{\Spec}{Spec}

\newcommand{\olot}{\operatorname{\ol\ot}}

\newcommand\braid{c}
\newcommand\mult{m}
\newcommand\multwo{\nabla}

\def\gr{\operatorname{gr}\nolimits}

\def\coh{\operatorname{H}\nolimits}
\def\Tot{\operatorname{Tot}\nolimits}

\DeclareMathOperator{\Span}{Span}

\newcommand{\HH}{\rm {H}}

\newcommand{\Z}{\mathbb Z}
\newcommand{\N}{\mathbb N}
\newcommand{\ot}{\otimes}

\newcommand{\wt}{\widetilde}
\newcommand{\bu}{\bullet}
\newcommand{\YD}{{}^{\Gamma}_{\Gamma}{\mathcal{YD}}}
\newcommand{\leer}{\text{---}}
\newcommand{\diag}{\operatorname{diag}}
\newcommand\ol{\overline}
\newcommand\inv{^{-1}}
\newcommand\ihom{\underline\hom}

\theoremstyle{plain}
\newtheorem{thm}{Theorem}[section]
\newtheorem{cor}[thm]{Corollary}
\newtheorem{lemma}[thm]{Lemma}
\newtheorem{lem}[thm]{Lemma}

\newtheorem{defn}[thm]{Definition}
\newtheorem{rem}[thm]{Remark}
\newtheorem{remark}[thm]{Remark}
\newtheorem{example}[thm]{Example}
\numberwithin{equation}{thm}

\title[Pointed Hopf algebras]
{Cohomology of finite dimensional pointed Hopf algebras}
\author{M.\ Mastnak}
\address{Department of Mathematics and Computer Science, Saint Mary's University,
Halifax, NS B3H 3C3, Canada}
\email{mmastnak@cs.smu.ca}
\author{J.\ Pevtsova}
\address{Department of Mathematics, University of Washington, Seattle, WA 98195, USA}
\email{julia@math.washington.edu}
\author{P.\ Schauenburg}
\address{Mathematisches Institut der Universit\"at M\"unchen, Theresienstr.\ 39, 80333 M\"unchen, Germany}
\email{schauenburg@math.lmu.de}
\author{S.\ Witherspoon}
\address{Department of Mathematics, Texas A\&M University, College Station, TX 77843, USA}
\email{sjw@math.tamu.edu}
\date{July 29, 2009}

\begin{document}

\maketitle

\begin{abstract}
We prove finite generation of the
cohomology ring of any finite dimensional pointed Hopf algebra,
having abelian group of grouplike elements,
under some mild restrictions on the group order.
The proof uses the recent classification by Andruskiewitsch and
Schneider of such Hopf algebras.
Examples include all of Lusztig's small quantum groups, whose
cohomology was first computed explicitly by Ginzburg and Kumar,
as well as many new pointed Hopf algebras.

We also show that in general the cohomology ring of a Hopf algebra in
a braided category is braided commutative.
As a consequence we obtain some further information about the structure
of the cohomology ring of a finite dimensional pointed Hopf
algebra and its related Nichols algebra.
\end{abstract}

\section{Introduction}

Golod \cite{Go}, Venkov \cite{V}, and Evens \cite{E}
proved that the cohomology ring of a finite
group with coefficients in a field of positive characteristic is finitely generated. This fundamental result opened the door
to using geometric methods in the study of cohomology and modular representations
of finite groups. This  geometric approach was pioneered by Quillen \cite{Q} and
expanded by Carlson \cite{C}, and Avrunin and Scott \cite{AvS}.
Friedlander and Suslin \cite{FS} vastly generalized the
result of Venkov and Evens, proving that the cohomology ring of any finite
group scheme (equivalently, finite dimensional cocommutative Hopf algebra) over a field of positive characteristic 
is finitely generated.
In a different direction, Ginzburg and Kumar \cite{GK} computed the cohomology
ring of each of Lusztig's small quantum groups $u_q({{\mathfrak g}})$, defined over $\mathbb C$,  under
some restrictions on the parameters; it is the coordinate ring of the nilpotent
cone of the Lie algebra $\mathfrak g$, and thus is finitely generated. The result is parallel to an
earlier result of Friedlander and  Parshall \cite{FPa}  who computed the cohomology
ring of a restricted Lie algebra in positive characteristic. 
Recently, Bendel, Nakano, Parshall, and Pillen \cite{BNPP}  calculated the cohomology  of
a small  quantum group under significantly loosened restrictions
on  the parameters. The broad common feature  of these works is that they investigate cohomology of certain 
nonsemisimple finite dimensional Hopf  algebras. Hence, these results lead one to ask whether the cohomology ring of any finite
dimensional Hopf algebra is finitely generated.
A positive answer
would simultaneously generalize the known results for
cocommutative Hopf algebras and for small quantum groups.
More generally, Etingof and Ostrik \cite{etingof-ostrik03} conjectured
finite generation of cohomology in the context of finite tensor categories.

In this paper,
we begin the task of proving this conjecture for more general classes
of noncocommutative Hopf algebras over a field of characteristic 0.
We prove finite generation of the cohomology
ring of any finite dimensional pointed Hopf algebra, with abelian group
of grouplike elements,
under some mild restrictions on the group order.
Pointed Hopf algebras are precisely those whose coradicals
are group algebras, and in turn these
groups determine a large part of their structure.
We use the recent classification of these Hopf algebras
by Andruskiewitsch and Schneider \cite{AS06}.
Each has a presentation by generators and relations similar to those
of the small quantum groups, yet they are much more general.
Due to this similarity, some of the techniques of Ginzburg and Kumar
\cite{GK} yield results in this general setting.
However some differences do arise,
notably that cohomology may no longer vanish in odd degrees,
and that useful connections to Lie algebras have not been developed.
Thus we must take a somewhat different approach in this paper,
which also yields new proofs of some of the results in \cite{GK}.
Since we have less information available in this general setting,
we prove finite generation without computing the full
structure of the cohomology ring.
However we do explicitly identify a subalgebra over which the cohomology
is finite, and establish some results about its structure.

These structure results follow in part from our general result
in Section \ref{bgcc}, that the cohomology ring of a Hopf algebra in
a braided category is always braided graded commutative.
This generalizes the well-known result
that the cohomology ring of a Hopf algebra is graded commutative,
one proof of which follows from the existence of two definitions of
its multiplication.
We generalize that proof, giving a braided version of the Eckmann-Hilton argument,
from which follows the braided graded commutativity result.
We  apply this to a Nichols algebra in Section \ref{graded},
thus obtaining some details about the structure of the finitely generated
cohomology ring of the corresponding pointed Hopf algebra.

Of course one hopes for  results in yet greater generality.
However the structure of finite dimensional noncommutative, noncocommutative, nonsemisimple
Hopf algebras, other than those treated in this paper, is largely unknown.
There are a very small number of known (nontrivial) finite dimensional pointed Hopf
algebras having nonabelian groups of grouplike elements (see for
example \cite{AF,AG}).
Even fewer examples are known of nonpointed, nonsemisimple Hopf algebras \cite{Men}.
To prove finite generation of cohomology in greater generality,
it may be necessary to find general techniques, such as the embedding
of any finite group scheme into GL$_n$ used by Friedlander and Suslin \cite{FS},
rather than depending  on structural knowledge of the Hopf algebras as
we do here.

Our proof of finite generation
is a two-step reduction to type $A_1\times \cdots\times A_1$,
in which case the corresponding Nichols algebra is a quantum complete
intersection.
For these algebras, we compute cohomology explicitly via a resolution constructed
in Section \ref{sec-step0}.
This resolution is adapted from \cite{BG,HK} where similar algebras were
considered (but in positive characteristic).
Each of our two reduction steps involves a spectral sequence associated to
an algebra filtration. In the first step (Section \ref{graded}) a Radford
biproduct of the form ${\mathcal B}(V)\# k\Gamma$, for a group $\Gamma$ and
Nichols algebra ${\mathcal B}(V)$, has a filtration for which
the associated graded algebra has type $A_1\times \cdots \times A_1$.
This filtration is generalized from De Concini and Kac \cite{DCK}.
We identify some permanent cycles and apply a lemma adapted from
Friedlander and Suslin \cite{FS} to conclude finite generation.
In the second reduction step (Section \ref{sec-step2}), any of
Andruskiewitsch and Schneider's pointed Hopf algebras $u({\mathcal D},\lambda,\mu)$
has a filtration for which the associated graded algebra is a Radford biproduct,
whose cohomology features in Section \ref{graded}.
Again we identify some permanent cycles and conclude finite generation.
As a corollary we show that the Hochschild cohomology ring of $u({\mathcal D},\lambda,\mu)$
is also finitely generated.

The first and last authors thank Ludwig-Maximilians-Universit\"at
M\"unchen for its hospitality during the first stages of this
project. The second author thanks MSRI for its hospitality and
support during the final stage of this work. The first author was
supported by an NSERC postdoctoral fellowship. The second author
was partially supported by the NSF grants DMS-0629156 and
DMS-0800940. The third author was supported by {\em Deutsche
Forschungsgemeinschaft} through a Heisenberg Fellowship. The last
author was partially supported by NSF grants DMS-0443476 and
DMS-0800832, and NSA grant H98230-07-1-0038. The last author
thanks D.\ J.\ Benson for very useful conversations.

\section{Definitions and Preliminary Results}\label{pointed}
Let $k$ be a field, usually assumed to be algebraically closed and
of characteristic 0.
All tensor products are over $k$ unless otherwise indicated.
Let $\Gamma$ be a finite group.

\subsection*{Hopf algebras in Yetter-Drinfeld categories}
A {\em Yetter-Drinfeld module} over $k\Gamma$ is a $\Gamma$-graded vector space
$V=\oplus_{g\in\Gamma}V_g$ that is also a $k\Gamma$-module for which
$g\cdot V_h = V_{ghg^{-1}}$ for all $g,h\in\Gamma$.
The grading by the group $\Gamma$ is equivalent to a $k\Gamma$-comodule
structure on $V$, that is a map $\delta: V\rightarrow k\Gamma \ot V$,
defined by $\delta(v)=g\ot v$ for all $v\in V_g$.
Let $\YD$ denote the category of all Yetter-Drinfeld modules over $k\Gamma$.
The category has a tensor product: $(V\ot W)_g = \oplus_{xy=g}V_x\ot W_y$ for all
$g\in\Gamma$, and $\Gamma$ acts diagonally on $V\ot W$, that is, $g(v\ot w)
=gv\ot gw$ for all $g\in \Gamma$, $v\in V$, and $w\in W$.
There is a {\em braiding} $c: V\ot W\stackrel{\sim}{\rightarrow} W\ot V$ for all $V,W\in \YD$
as follows:
Let $g\in\Gamma$, $v\in V_g$, and $w\in W$.
Then
$
   c(v\ot w) = gw\ot v.
$
Thus $\YD$ is a braided monoidal category.
(For details on the category $\YD$, including the connection to
Hopf algebras recalled below, see for example \cite{AS06}.)

Let $\mathcal C$ be any braided monoidal category. For simplicity, we
will always assume the tensor product is strictly associative
with strict unit object $I$. An {\em algebra in} $\mathcal C$ is
an object $R$ together with morphisms $u: I\rightarrow R$ and $m:
R\ot R\rightarrow R$ in $\mathcal C$, such that $m$ is associative
in the sense that $m(m\ot 1_R)=m(1_R\ot m)$, and $u$ is a unit in
the sense that $m(u\ot 1_R)=1_R=m(1_R\ot u)$. The definition of a
coalgebra in $\mathcal C$ is similar, with the arrows going in the
opposite direction. Thus, an algebra (resp.\ coalgebra) in $\YD$
is simply an ordinary algebra (resp.\ coalgebra) with
multiplication (resp.\ comultiplication) a graded and equivariant
map. A {\em braided Hopf algebra in} $\mathcal C$ is an algebra
as well as coalgebra in $\mathcal C$ such that its
comultiplication and counit are algebra morphisms, and such that
the identity morphism $\id\colon R\to R$ has a convolution inverse
$s$ in $\mathcal C$. When we say that the comultiplication
$\Delta\colon R\to R\ot R$ should be an algebra morphism, the
braiding $c$ of $\mathcal C$ arises in the definition of the
algebra structure of $R\ot R$ (so in particular, a braided Hopf
algebra in $\YD$ is not an ordinary Hopf algebra). More generally,
if $A,B$ are two algebras in $\mathcal C$, their tensor product
$A\ot B$ is defined to have multiplication $m_{A\ot B}=(m_A\ot
m_B)(1_A\ot c\ot 1_B)$.

An example of a braided Hopf algebra in $\YD$ is the {\em Nichols algebra}
${\mathcal B}(V)$ associated to a Yetter-Drinfeld module $V$ over $k \Gamma$;
${\mathcal B}(V)$ is the quotient of the tensor algebra $T(V)$ by the
largest homogeneous braided Hopf ideal generated by homogeneous elements
of degree at least 2. For details, see \cite{AS02,AS06}.
In this paper, we need only the structure of ${\mathcal B}(V)$ in some
cases as are explicitly described below.

If $R$ is a braided Hopf algebra in $\YD$, then its {\em Radford
biproduct} (or {\em bosonization}) $R\# k\Gamma$ is a Hopf algebra
in the usual sense (that is, a Hopf algebra in ${}^k_k{\mathcal
{YD}}$). As an algebra, $R\# k\Gamma$ is just a {\em skew group
algebra}, that is, $R\# k\Gamma$ is a free $R$-module with basis
$\Gamma$ on which multiplication is defined by $
   (rg)(sh) = r(g\cdot s) gh
$
for all $r,s\in R$ and $g,h\in G$.
Comultiplication is given by
$
  \Delta(rg)=r^{(1)} (r^{(2)})_{(-1)} g \ot (r^{(2)})_{(0)} g ,
$
for all $r\in R$ and $g\in \Gamma$, where $\Delta(r)=\sum r^{(1)}\ot r^{(2)}$ in $R$
as a Hopf algebra in $\YD$ and $\delta(r)=\sum r_{(-1)}\ot r_{(0)}$ denotes the
$k\Gamma$-comodule structure.

\subsection*{The pointed Hopf algebras of Andruskiewitsch and Schneider}
A {\em pointed} Hopf algebra $H$ is one for which all simple comodules are
one dimensional.
This is equivalent to the condition $H_0=k\Gamma$ where $\Gamma = G(H)$ is the
group of grouplike elements of $H$ and $H_0$ is the coradical of $H$
(the initial term in the coradical filtration).

The Hopf algebras of Andruskiewitsch and Schneider in \cite{AS06} are pointed, and
are deformations of Radford biproducts.
They depend on the following data:
Let $\theta$ be a positive integer.
Let $(a_{ij})_{1\leq i,j\leq \theta}$ be a
{\em Cartan matrix of finite type}, that is the Dynkin
diagram of $(a_{ij})$ is a disjoint union of copies of some of the diagrams
$A_{\bullet}$, $B_{\bullet}$, $C_{\bullet}$, $D_{\bullet}$, $E_6$, $E_7$, $E_8$, $F_4$, $G_2$.
In particular,
$a_{ii}=2$ for $1\leq i\leq \theta$,
$a_{ij}$ is a nonpositive integer for $i\neq j$, and $a_{ij}=0$ implies
$a_{ji}=0$.
Its {\em Dynkin diagram} is a graph with vertices labeled $1,\ldots,\theta$:
The vertices $i$ and $j$ are connected by
$a_{ij}a_{ji}$ edges, and if $|a_{ij}|>|a_{ji}|$, there is an arrow pointing
from $j$ to $i$.

Now assume $\Gamma$ is {\em abelian},  and denote by
$\hat{\Gamma}$ its dual group of characters.
For each $i$, $1\leq i\leq \theta$, choose $g_i\in \Gamma$ and
$\chi_i\in\hat{\Gamma}$
such that $\chi_i(g_i)\neq 1$  and
\begin{equation}
\label{cartan}
   \chi_j(g_i)\chi_i(g_j)=\chi_i(g_i)^{a_{ij}}
\end{equation}
(the {\em Cartan condition}) holds for $1\leq i,j\leq \theta$.
Letting $q_{ij}=\chi_j(g_i)$, this becomes $q_{ij}q_{ji}=q_{ii}^{a_{ij}}$.
Call
\begin{equation}\label{data}
  {\mathcal D} = (\Gamma, (g_i)_{1\leq i\leq \theta},
   (\chi_i)_{1\leq i\leq \theta}, (a_{ij})_{1\leq i,j\leq \theta})
\end{equation}
a {\em datum of finite Cartan type} associated to $\Gamma$ and $(a_{ij})$.
The Hopf algebras of interest will be generated as algebras by $\Gamma$
and symbols $x_1,\ldots,x_{\theta}$.

Let $V$ be the vector space with basis $x_1,\ldots,x_{\theta}$.
Then $V$ has a structure of a Yetter-Drinfeld module over
$k\Gamma$: $V_g=\Span_k\{x_i\mid g_i=g\}$ and
$g(x_i)=\chi_i(g)x_i$ for $1\leq i\leq \theta$ and $g\in\Gamma$.
This induces the structure of an algebra in $\YD$ on the tensor
algebra $T(V)$. In particular $\Gamma$ acts by automorphisms on
$T(V)$, and
$T(V)$ is a $\Gamma$-graded algebra in which $x_{i_1}\cdots
x_{i_s}$ has degree $g_{i_1}\cdots g_{i_s}$.
The braiding $c\colon T(V)\ot
T(V) \rightarrow T(V)\ot T(V)$ is induced by $c(x_i\ot y) =
g_i(y)\ot x_i$. Moreover, $T(V)$ can be made a braided Hopf
algebra in $\YD$ if we define comultiplication as the unique
algebra map $\Delta\colon T(V)\to T(V)\ot T(V)$ satisfying
$\Delta(v)=v\ot 1+1\ot v$ for all $v\in V$. We define the braided
commutators
$$
  \ad_c(x_i)(y) = [x_i,y]_c := x_i y - g_i(y) x_i,
$$
for all $y\in T(V)$, and similarly in quotients of $T(V)$
by homogeneous ideals.

Let $\Phi$ denote the root system corresponding to $(a_{ij})$, and let $\Pi$
denote a fixed set of simple roots.
If $\alpha_i, \alpha_j\in \Pi$, write $i\sim j$ if the corresponding
vertices in the Dynkin diagram of $\Phi$ are in the same connected component.
Choose scalars $\lambda = (\lambda_{ij})_{1\leq i<j\leq \theta, \ i\not\sim j}$,
called {\em linking parameters}, such that
\begin{equation*}
   \lambda_{ij}=0 \ \ \mbox{ if } \ g_ig_j=1 \ \mbox{ or } \ \chi_i\chi_j\neq \varepsilon,
\end{equation*}
where $\varepsilon$ is the identity element in the dual group $\hat{\Gamma}$
(equivalently $\varepsilon$ is the counit on $k\Gamma$, $\varepsilon(g)=1$ for all
$g\in \Gamma$).
The Hopf algebra $U({\mathcal D},\lambda)$ defined by Andruskiewitsch and
Schneider \cite{AS06} is the quotient of $T(V)\# k\Gamma$ by relations
corresponding to the equations
\begin{eqnarray*}
  \mbox{ {\em (group action)}}\hspace{1.19in} \ \
  gx_ig^{-1} &=& \chi_i(g) x_i  \ \ \ \ \ (g\in\Gamma, 1\leq i\leq \theta),
\\
  \mbox{ {\em (Serre relations)}} \hspace{.3in}\ \
 (\ad_c(x_i))^{1-a_{ij}} (x_j) &=& 0 \ \ \ \ \  (i\neq j, \ i\sim j),
\\
  \mbox{{\em (linking relations)}} \hspace{.51in}\ \
 (\ad_c(x_i))(x_j)& =& \lambda_{ij} (1-g_ig_j)  \ \ \ \ \
    (i<j, \ i\not\sim j).
\end{eqnarray*}
The coalgebra structure of $U({\mathcal{D}},\lambda)$ is given by
$$
 \Delta(g) = g\ot g , \ \ \ \Delta(x_i)=x_i\ot 1 + g_i\ot x_i,
$$
$\varepsilon(g)=1$, $\varepsilon(x_i)=0$, $s(g)=g^{-1}$, $s(x_i)=-g_i^{-1}x_i$,
for all $g\in \Gamma$, $1\leq i\leq \theta$.

\begin{example}
\label{uqg} {\em Let $\mathfrak g$ be a simple Lie algebra of rank
$n$, let $(a_{ij})$ be {\em two} block diagonal copies of the
corresponding Cartan matrix, and let $\theta=2n$. Let  $q$ be a
primitive $\ell$th root of unity, $\ell$ odd. Let $\Gamma
=(\Z/\ell\Z)^n$ with generators $g_1,\ldots,g_n$ and define
$\chi_i\in \hat{\Gamma}$ by $\chi_i(g_j) = q^{\langle \alpha_i, \
\alpha_j\rangle}$. Let $g_{i+n}=g_i$ and $\chi_{i+n}=\chi_i^{-1}$
for $1\leq i\leq n$. Let $\lambda_{ij} = (q^{-1}-q)^{-1}
\delta_{j, i+n}$ for $1\leq i<j\leq n$. Then $U({\mathcal
D},\lambda)$ is a quotient of the quantum group $U_q({\mathfrak
g})$. An epimorphism $U_q(\mathfrak g)\to U(\mathcal D,\lambda)$
is given by $K_i\mapsto g_i$, $E_i\mapsto x_i$ and $F_i\mapsto
x_{i+n}g_i\inv $.}
\end{example}

The Hopf algebra $U({\mathcal D},\lambda)$ has finite dimensional
quotients that we discuss next. As in \cite{AS06} we make the
assumptions:
\begin{equation}
\label{assumptions}
\begin{array}{l}
  \mbox{{\em The order of $\chi_i(g_i)$ is odd for all $i$,} }\\
  \mbox{{\em and is prime to 3 for all $i$ in a connected component of type $G_2$.}}
\end{array}
\end{equation}
It then follows from the Cartan condition (\ref{cartan}) that
the order of $\chi_i(g_i)$ is the same as the order of $\chi_j(g_j)$ if $i\sim j$.
That is, this order is constant in each connected component
$J$ of the Dynkin diagram;
denote this common order by $N_J$.
It will also be convenient to denote it by $N_{\beta_j}$ or $N_j$
for each positive root  $\beta_j$ in $J$ (this standard notation is
defined below).
Let $\alpha\in\Phi^+$, $\alpha=\sum_{i=1}^{\theta} n_i\alpha_i$, and
let $\Ht(\alpha)=\sum_{i=1}^{\theta}n_i$,
\begin{equation}\label{gbeta-chibeta}
 g_{\alpha}=\prod_{i=1}^{\theta} g_i^{n_i}, \ \ \mbox{ and }
 \ \ \chi_{\alpha} =\prod_{i=1}^{\theta} \chi_i^{n_i}.
\end{equation}
There is a unique connected component $J_{\alpha}$ of the Dynkin
diagram of $\Phi$ for which $n_i\neq 0$ implies $i\in J_{\alpha}$.
We write $J=J_{\alpha}$ when it is clear which $\alpha$
is intended.

Let $W$ be the  Weyl group of the root system $\Phi$.
Let $w_0=s_{i_1}\cdots s_{i_r}$
be a reduced decomposition of the longest
element $w_0\in W$ as a product of simple reflections.
Let
$$
  \beta_1=\alpha_{i_1}, \ \ \beta_2=s_{i_1}(\alpha_{i_2}),\ \ldots, \
 \beta_r=s_{i_1} s_{i_2}\cdots s_{i_{r-1}} (\alpha_{i_r}).
$$
Then $\beta_1,\ldots,\beta_r$ are precisely the positive roots $\Phi^+$
\cite{L}.
Corresponding root vectors $x_{\beta_j}\in U({\mathcal D},\lambda)$
are defined in the same way as for the traditional quantum groups:
In case $\mathcal D$ corresponds to the data for a quantum group
$U_q({\mathfrak{g}})$ (see Example \ref{uqg}), let
$$
   x_{\beta_j}=T_{i_1}T_{i_2}\cdots T_{i_{j-1}}(x_{i_j}),
$$
where the $T_{i_j}$ are Lusztig's algebra automorphisms of
$U_q({\mathfrak{g}})$ \cite{L}.
The $x_{\beta_j}$ may be expressed as iterated braided commutators.
If $\beta_j$ is a simple root $\alpha_l$, then $x_{\beta_j}=x_l$.
In our more general setting, as in \cite{AS06},
define the $x_{\beta_j}$ to be the analogous iterated
braided commutators.
More precisely, we describe how to obtain these elements in case the
Dynkin diagram is connected; in the general case construct the
root vectors separately for each connected component.
Let $I$ be the ideal of $T(V)$ generated by elements corresponding
to the Serre relations. Then
\cite[Lemmas 1.2 and 1.7]{AS06} may be used to obtain a
{\em linear} isomorphism between $T(V)/I$ and the upper triangular
part of some $U_q({\mathfrak g})$ with the same Dynkin diagram.
This linear isomorphism preserves products and braided commutators
up to nonzero scalar multiples, and thus yields $x_{\beta_j}$ in our general
setting as an iterated braided commutator.
(See the proof of
Lemma \ref{filtration} below for more details on this linear isomorphism.)

Choose scalars $(\mu_{\alpha})_{\alpha\in\Phi^+}$, called {\em root
vector parameters}, such that
\begin{equation}
\label{mualpha-zero}
  \mu_{\alpha}=0 \ \mbox{ if } \ g_{\alpha}^{N_{\alpha}}=1 \
  \mbox{ or } \ \chi_{\alpha}^{N_{\alpha}}\neq \varepsilon.
\end{equation}
The finite dimensional Hopf algebra $u({\mathcal D},\lambda,\mu)$ is the
quotient of $U({\mathcal D},\lambda)$ by the ideal generated by
all
\begin{equation}
\label{root-relns}
  \hspace{-.5cm}
 \mbox{ {\em (root vector relations)}} \hspace{1.8cm}
   x_{\alpha}^{N_{\alpha}} - u_{\alpha}(\mu) \ \ \ \ \ \ (\alpha\in \Phi^+)
\end{equation}
where $u_{\alpha}(\mu)\in k\Gamma$ is defined inductively on $\Phi^+$
in \cite[Defn.\ 2.14]{AS06}.
In this paper, we do not need the details of the construction
of the elements $u_{\alpha}(\mu)$ in the group algebra. We only need
the fact that if $\mu_{\alpha}=0$ for all $\alpha\in\Phi^+$, then
$u_{\alpha}(\mu)=0$ for all $\alpha\in\Phi^+$ \cite{AS06} (see for example
Lemma \ref{mu-zero} below).  It is interesting to note that if
 $\alpha$ is a {\em simple} root, then $u_{\alpha}(\mu)
:=\mu_{\alpha}(1-g_{\alpha}^{N_{\alpha}})$.

\begin{example}
{\em
Let $\mathcal D$, $\lambda$ be the data from Example \ref{uqg}.
Then there is an isomorphism $u_q({\mathfrak g})\simeq u({\mathcal D},\lambda,0)$,
induced by the epimorphism $U_q({\mathfrak g})\rightarrow U({\mathcal
D},\lambda)$ given in that example.}
\end{example}

The following theorem is \cite[Classification Theorem 0.1]{AS06}, and
requires $k$ to be algebraically closed of characteristic 0.

\begin{thm}[Andruskiewitsch-Schneider]
The Hopf algebras $u({\mathcal{D}}, \lambda,\mu)$ are finite dimensional
and pointed. Conversely, if $H$ is a finite dimensional pointed Hopf algebra having abelian
group of grouplike elements with order not divisible by primes less than 11,
then $H\simeq u({\mathcal{D}},\lambda,\mu)$ for some ${\mathcal D}$, $\lambda$,
$\mu$.
\end{thm}

\subsection*{Two filtrations}
Note that $u({\mathcal D},\lambda,\mu)$ is a (coradically) filtered Hopf algebra,
with $\deg (x_i) = 1$ ($1\leq i\leq \theta$) and $\deg (g) = 0$ ($g\in \Gamma$).
The associated graded Hopf algebra is isomorphic to the graded Hopf algebra
$u({\mathcal {D}},0,0)$.
There is an isomorphism $u({\mathcal {D}},0,0)\simeq {\mathcal{B}}(V)\# k\Gamma$,
the Radford biproduct (or bosonization) of the Nichols algebra
${\mathcal{B}}(V)$ of the Yetter-Drinfeld module $V$ over $k\Gamma$.
For details, see \cite{AS06}.
Note that ${\mathcal B}(V)$ is isomorphic to the subalgebra of $u({\mathcal {D}},0,0)$
generated by $x_1,\ldots,x_{\theta}$.

In Section \ref{sec-step2}, we prove that the cohomology of $u({\mathcal D},\lambda,\mu)$
is finitely generated by using a spectral sequence relating its cohomology to that of
${\mathcal B}(V)\# k\Gamma$.
In Section~\ref{graded}, we prove
that the cohomology of this (coradically) graded bialgebra ${\mathcal{B}}(V)\#
k\Gamma$ is finitely generated, by using a spectral sequence relating
its cohomology to that of a much simpler algebra: We put a
different filtration on it (see Lemma \ref{filtration} below) for which the
associated graded algebra has type $A_1\times\cdots\times A_1$.
(That is, the Dynkin diagram is a disjoint union of diagrams $A_1$.)
In Section \ref{sec-step0}, we give the cohomology for type $A_1\times\cdots\times A_1$
explicitly as a special case of the cohomology of a quantum complete intersection.

By \cite[Thm.\ 2.6]{AS06}, the Nichols algebra ${\mathcal B}(V)$ has PBW basis all
\begin{equation}
\label{PBW}
  {\bf x}^{\bf a}=x_{\beta_1}^{a_1}\cdots x_{\beta_r}^{a_r} \ \ \
  (0\leq a_i < N_{\beta_i}),
\end{equation}
and further
\begin{equation}
\label{Ncentral}
[x_{\alpha},x_{\beta}^{N_{\beta}}]_c = 0
\end{equation}
for all $\alpha,\beta\in\Phi^+$.
As in \cite{DCK}, put a total order on the PBW basis elements
as follows: The degree of such an element is
$$
  d(x_{\beta_1}^{a_1}\cdots x_{\beta_r}^{a_r})
   =\left(\prod a_i\Ht(\beta_i), a_r,\ldots,a_1\right)\in \N^{r+1}.
$$
Order the elements (\ref{PBW}) lexicographically by degree where
$$
  (0,\ldots,0,1)<(0,\ldots,0 ,1,0)<\cdots <(1,0,\ldots ,0).
$$

\begin{lem}\label{filtration}
In the Nichols algebra ${\mathcal B}(V)$,
for all $i<j$,
$$
   [x_{\beta_j},x_{\beta_i}]_c = \sum_{{\bf a}\in\N^p}
    \rho_{\bf a}{\bf x}^{\bf a}
$$
where the $\rho_{\bf a}$ are scalars for which $\rho_{\bf a}=0$
unless $d({\bf x}^{\bf a}) < d(x_{\beta_i}x_{\beta_j})$.
\end{lem}

\begin{proof}
First note that if $i\not\sim j$, then $[x_{\beta_j},x_{\beta_i}]_c=0$
by the Serre relations, since $a_{ij}=0$.
Thus we may assume now that $i\sim j$.
Then the lemma is just the translation of \cite[Lemma 1.7]{DCK}, via twisting
by a cocycle and specializing $q$, into this more general setting.
The twisting method is described in \cite{AS06} and is used there
to prove that (\ref{PBW}) is a basis of ${\mathcal B}(V)$.
In particular,  \cite[Lemma 2.3]{AS06} states that there exist
integers $d_i\in\{1,2,3\}$, $1\leq i\leq \theta$, and $q\in k$ such that
for all $1\leq i,j\leq \theta$, $q_{ii}=q^{2d_i}$ and $d_ia_{ij}=d_ja_{ji}$.
This allows us to define a matrix
$(q_{ij}')$ by $q_{ij}'=q^{d_ia_{ij}}$, so that
$q_{ij}q_{ji}=q_{ij}'q_{ji}'$ and $q_{ii}=q_{ii}'$.
Let $V'\in\YD$ with basis $x_1'\ldots,x_{\theta}'$ and structure given
by the data $\{q_{ij}'\}$.
Now $\{q_{ij}'\}$ are classical quantum group parameters for the positive
part of $U_q({\mathfrak g})$, where $\mathfrak g$ corresponds to the Dynkin diagram.
Thus \cite[Lemma 1.7]{DCK} applies. Further, \cite[Lemma 1.2]{AS06} states
that there is a cocycle $\sigma: \Z\Pi\times\Z\Pi\rightarrow k^{\times}$
and a $k$-linear isomorphism
$\phi: T(V)\rightarrow T(V')$ with $\phi(x_i)=x_i'$ and
such that if $x\in T(V)_{g_i}$ and $y\in T(V)_{g_j}$, then
$$
  \phi(xy) = \sigma(\alpha_i,\alpha_j) \phi(x)\phi(y) \ \ \ \mbox{ and } \ \ \
  \phi([x,y]_c) = \sigma(\alpha_i,\alpha_j)[\phi(x),\phi(y)]_{c'}.
$$
Thus $\phi$  preserves the PBW basis up to nonzero scalar multiples,
and it preserves the total order on the PBW basis.
Since it also preserves braided commutators up to nonzero scalar,
the lemma now holds as a consequence of the same result
\cite[Lemma 1.7]{DCK} for the parameters $\{q_{ij}'\}$.
\end{proof}

By Lemma \ref{filtration},
the above ordering induces a filtration $F$ on ${\mathcal B}(V)$
for which the associated graded algebra
$\Gr {\mathcal B}(V)$ has relations $[x_{\beta_j},x_{\beta_i}]_c =0$ for all $i<j$,
and $x_{\beta_i}^{N_i} =0$.
In particular, $\Gr {\mathcal B}(V)$ is of type $A_1\times\cdots\times A_1$.
We may put a corresponding Hopf algebra structure on $(\Gr {\mathcal B}(V))\# k\Gamma$
as follows.
If $\beta_i = \sum_{j=1}^{\theta} n_j\alpha_j$, let $g_{\beta_i}=g_1^{n_1}\cdots
g_{\theta}^{n_{\theta}}$ as in (\ref{gbeta-chibeta}).
Now identify $\beta_1,\ldots,\beta_r$ with the simple roots of type
$A_1\times\cdots\times A_1$.
For $i<j$, define $q_{\beta_i\beta_j} = \chi_{\beta_j}(g_{\beta_i})$,
$q_{\beta_j\beta_i} = q_{\beta_i\beta_j}^{-1}$, and $q_{\beta_i\beta_i}=1$.
Then the Cartan
condition (\ref{cartan}) holds for the scalars $q_{\beta_i\beta_j}$ in
this type $A_1\times\cdots\times A_1$, and
$(\Gr {\mathcal B}(V))\# k\Gamma$ is a Hopf algebra for which
$$\Delta(x_{\beta_i}) = x_{\beta_i}\ot 1 + g_{\beta_i}\ot x_{\beta_i}.
$$

For our spectral sequence constructions, we rewrite
the total order on the PBW basis elements (\ref{PBW})
explicitly as an indexing by
positive integers.
We may set
$$
  \deg(x_{\beta_1}^{a_1}\cdots x_{\beta_r}^{a_r})=  N_{\beta_1}\cdots N_{\beta_r}\prod a_i\Ht (\beta_i) +
   N_{\beta_1}\cdots N_{\beta_{r-1}} a_r + \cdots + N_{\beta_1} a_2 + a_1.
$$
A case-by-case argument shows that  $d(x_{\beta_1}^{a_1}\cdots x_{\beta_r}^{a_r})
<d(x_{\beta_1}^{b_1}\cdots x_{\beta_r}^{b_r})$ if, and only if,
$\deg(x_{\beta_1}^{a_1}\cdots x_{\beta_r}^{a_r})<\deg(x_{\beta_1}^{b_1}\cdots
x_{\beta_r}^{b_r})$.

\subsection*{Hochschild cohomology}
In this paper, we are interested in the cohomology ring $\coh^*(u({\mathcal D},\lambda,\mu),k)
:=\Ext^*_{u({\mathcal D},\lambda,\mu)}(k,k)$, where $k$ is the trivial module given
by the counit $\varepsilon$.
If $A$ is any $k$-algebra with an augmentation $\varepsilon: A\rightarrow k$,
note that $\Ext^*_A(k,k)$ is isomorphic to Hochschild
cohomology with trivial coefficients, $\Ext^*_{A\ot A^{op}}(A,k)$:
This is due to an equivalence of the bar complexes for computing these Ext algebras.
That is, these bar complexes  are each equivalent to
the reduced complex
\begin{equation}\label{cdot}
  C^{\bullet}: \hspace{.5in}
   0\rightarrow\Hom_k(k,k)\rightarrow \Hom_k(A^+,k)\rightarrow\Hom_k((A^+)^{2},k)\rightarrow
\cdots,
\end{equation}
where $A^+=\ker \varepsilon$ is the augmentation ideal of $A$.
The differential is given by
$\delta_{n+1}(f) = \sum_{i=0}^{n-1}(-1)^{i+1} f\circ (1^i\ot m\ot 1^{n-i-1})$
for all $f: (A^+)^{\ot n}\rightarrow k$.
We will exploit this equivalence.
This complex arises, for example, by applying $\Hom_A( - ,k)$ to the free $A$-resolution
of $k$:
$$
  \cdots\rightarrow A\ot (A^+)^{\ot 2} \stackrel{\partial_2}{\longrightarrow}
   A\ot A^+ \stackrel{\partial_1}{\longrightarrow} A\stackrel{\varepsilon}{\longrightarrow} k
  \rightarrow 0,
$$
where
\begin{equation}\label{free-res}
  \partial_i(a_0\ot\cdots\ot a_i)=\sum_{j=0}^{i-1} (-1)^j a_0\ot\cdots\ot a_j a_{j+1}
  \ot\cdots \ot a_i.
\end{equation}
Then $\partial^*_{n+1} =\delta_{n+1}$.
Equivalently, we may apply $\Hom_{A\ot A^{op}}( - , k)$ to the
free $A\ot A^{op}$-resolution of $A$:
$$
 \cdots\rightarrow A\ot (A^+)^{\ot 2}\ot A \stackrel{d_2}{\longrightarrow}
   A\ot A^+\ot A \stackrel{d_1}{\longrightarrow} A\ot A\stackrel{\varepsilon}{\longrightarrow} A
  \rightarrow 0,
$$
$d_i(a_0\ot \cdots\ot a_{i+1}) =\sum_{j=0}^i (-1)^j a_0\ot\cdots
  \ot a_ja_{j+1}\ot \cdots\ot a_{i+1}.
$

\subsection*{A finite generation lemma}
In Sections \ref{graded} and \ref{sec-step2}, we
will need the following general lemma adapted from \cite[Lemma 1.6]{FS}.
Recall that an element $a\in E_r^{p,q}$ is called a  {\it permanent cycle}
if $d_i(a) = 0$ for all $i \geq r$.

\begin{lemma}
\label{fingen}
{\rm (a)} Let $\xymatrix{E_{1}^{p,q} \ar@{=>}[r] &  E^{p+q}_\infty}$
be a multiplicative spectral sequence  of $k$-algebras concentrated
in the half plane $p+q \geq 0$,
and let $A^{*,*}$ be a bigraded commutative $k$-algebra concentrated in even
(total) degrees.
Assume that there exists a bigraded map of algebras
$\xymatrix{\phi: A^{*,*} \ar[r] & E_1^{*,*}}$ such that

\begin{enumerate}
\item[(i)] $\phi$ makes $E_1^{*,*}$ into a Noetherian $A^{*,*}$-module, and
\item[(ii)] the image of $A^{*,*}$ in $E_1^{*,*}$ consists  of  permanent cycles.
\end{enumerate}
Then $E^*_\infty$ is a Noetherian module over $\Tot(A^{*,*})$.

\vspace{0.1in}
\noindent
{\rm (b)} Let $\xymatrix{\wt E_1^{p,q} \ar@{=>}[r] &  \wt E^{p+q}_\infty}$ be a spectral
sequence that  is a bigraded module over the spectral sequence $E^{*,*}$.
Assume that $\wt E_1^{*,*}$ is a Noetherian module  over
$A^{*,*}$ where $A^{*,*}$ acts on $\wt E_1^{*,*}$ via the map $\phi$.
Then $\wt E^*_\infty$ is a finitely generated $E^*_\infty$-module.
\end{lemma}

\begin{proof}
Let $\Lambda^{*,*}_r \subset E_r^{*,*}$ be the bigraded subalgebra of permanent
cycles  in $E^{*,*}_r$.
Observe that $d_r(E_r^{*,*})$ is an $A^{*,*}$-invariant left ideal of $\Lambda^{*,*}_r$.
Indeed, let $ a\in A^{p,q}$ and  $x \in E^{s,t}_r$.
We have  $d_r(\overline{\phi(a)}x) = d_r(\overline{\phi(a)}) x +
(-1)^{p+q}\overline{\phi(a)} d_r(x) =
 \overline{\phi(a)} d_r(x)$ since $\overline{\phi(a)}\in A^{*,*}$ is assumed to be a
permanent  cycle of even total degree.  A similar computation  shows  that
$\Lambda_1^{*,*}$  is an $A^{*,*}$-submodule of $E_1^{*,*}$.  By induction,
 $\Lambda^{*,*}_{r+1} = \Lambda^{*,*}_r/d_r(E_r^{*,*})$ is an $A^{*,*}$-module
for any $r \geq 1$.
  We get a sequence  of surjective maps of $A^{*,*}$-modules:
\begin{equation}
\label{surj}
\xymatrix{ \Lambda^{*,*}_1 \ar@{->>}[r] &\ldots   \ar@{->>}[r] &
\Lambda^{*,*}_r \ar@{->>}[r] & \Lambda^{*,*}_{r+1} \ar@{->>}[r] & \ldots}
\end{equation}
Since $\Lambda^{*,*}_1$ is an $A^{*,*}$-submodule  of $E_1^{*,*}$,
it is Noetherian as an $A^{*,*}$-module. Therefore, the kernels of
the maps $\xymatrix{ \Lambda^{*,*}_1 \ar@{->>}[r] &
\Lambda^{*,*}_r}$ are Noetherian  for all $r\geq 1$. These kernels
form an increasing chain of submodules of $\Lambda^{*,*}_1$,
hence, by the Noetherian property, they  stabilize after finitely
many steps. Consequently, the
sequence~(\ref{surj}) stabilizes after finitely many steps. We
conclude that $\Lambda^{*,*}_r = E_\infty^{*,*}$ for some $r$.
Therefore $E_\infty^{*,*}$ is a Noetherian $A^{*,*}$-module. By
\cite[Proposition 2.1]{E}, $E^*_\infty$ is a Noetherian module
over $\Tot(A^{*,*})$, finishing the proof of (a).

Let $\wt\Lambda^{*,*}_r \subset \wt E_r^{*,*}$ be the subspace  of permanent cycles.
Arguing as above we can
show that $\wt\Lambda^{*,*}_r$ is an $A^{*,*}$-submodule  of $\wt E^{*,*}_r$, and,
moreover, that there exists $r$ such that
$\wt\Lambda^{*,*}_r = \wt E^{*,*}_\infty$. Hence, $\wt E^{*,*}_\infty$ is Noetherian
over $A^{*,*}$.  Applying  \cite[Proposition 2.1]{E}
once again, we conclude that $\wt E^*_\infty$ is Noetherian and hence
finitely generated over $\Tot(A^{*,*})$.   Therefore, it is finitely generated
over $E^*_\infty$.
\end{proof}

\section{Braided graded commutativity of cohomology}
\label{bgcc}

We will show that Hochschild cohomology of a braided Hopf algebra
in $\YD$ is a braided graded commutative algebra (more precisely,
its opposite algebra is). This will be done in a more general
categorical setting, so it will make sense to relax the overall
requirement that $k$ is a field.

To prepare for the main result of this section, we need some
information on the Alexander-Whitney map. We follow Weibel
\cite[Section 8.5]{W}. The case of a tensor product of simplicial
modules over a commutative ring (rather than simplicial objects in
a tensor category, or bisimplicial objects) can also be found in
Mac Lane \cite[Ch.~VII, Section~8]{maclane}.

Consider a bisimplicial object
$A=(A_{k\ell},\partial^h,\partial^v)$ in an abelian category.
Associated to it we have the diagonal simplicial object $\diag A$
with components $A_{nn}$; to this in turn we associate the usual
chain complex, whose differentials we denote by $d^{\diag}_n$.
Also associated to $A$ is the double complex $CA$, whose
horizontal and vertical differentials we denote by
$d^h_{k,\ell}\colon A_{k\ell}\to A_{k-1,\ell}$ and
$d^v_{k,\ell}\colon A_{k\ell}\to A_{k,\ell-1}.$ There is a natural
chain morphism, unique up to natural chain homotopy,
$$f\colon \diag A\to \Tot CA$$
for which $f_0$ is the identity. We note that $f$ being a chain
morphism means that the components $f_{k\ell}\colon A_{nn}\to
A_{k\ell}$ for $n=k+\ell$ satisfy
\[f_{k\ell}\, d_{n+1}^{\diag}=d_{k+1,\ell}^hf_{k+1,\ell}+d_{k,\ell+1}^vf_{k,\ell+1}.\]
One possible choice for $f$ is the Alexander-Whitney map, whose
components are
\[f_{k\ell}=\partial^h_{k+1}\dots\partial^h_n\underbrace{\partial^v_0\dots\partial^v_0}_k\colon A_{nn}\to A_{k\ell}\]
for $n=k+\ell$. We claim that
\[f'_{k\ell}=(-1)^{k\ell}\partial^h_0\dots\partial^h_0\partial^v_{\ell+1}\dots\partial^v_n\]
is another choice, and therefore chain homotopic to $f$. For this
it suffices to check that $f'$ (which is clearly natural) is a
chain morphism. Instead of doing this directly (by more or less
the same calculations usually done for the Alexander-Whitney map),
we consider the front-to-back dual version of $A$. We denote this
by $\tilde A$.
The horizontal and vertical face operators are $\tilde{\partial}^h_k
=\tilde{\partial}^h_{n-k}$ and $\tilde{\partial}^v_k
=\tilde{\partial}^v_{n-k}$.
The horizontal and vertical differentials in
$C\tilde A$ are $\tilde d^h_{k,\ell}=(-1)^kd^h_{k,\ell}$ and
$\tilde d^v_{k,\ell}=(-1)^\ell d^v_{k,\ell}$, the differentials in
$\diag\tilde A$ are $\tilde d^{\diag}_n=(-1)^n d^{\diag}_n$. The
Alexander-Whitney map $\tilde f$ of $\tilde A$ is related to $f'$
by $f'_{k\ell}=(-1)^{k\ell}\tilde f_{k\ell}$. Therefore
\begin{align*}
  f'_{k\ell}\, d^{\diag}_{n+1}&=(-1)^{k\ell}\tilde
  f_{k\ell}(-1)^{n+1}\tilde d^{\diag}_{n+1}\\
  &=(-1)^{k\ell+k+\ell+1}(\tilde d^h_{k+1,\ell}\tilde f_{k+1,\ell}+\tilde d^v_{k,\ell+1}\tilde
  f_{k,\ell+1})\\
  &=(-1)^{k+1}\tilde d^h_{k+1,\ell}(-1)^{(k+1)\ell}\tilde f_{k+1,\ell}+
  (-1)^{\ell+1}\tilde d^v_{k,\ell+1}(-1)^{k(\ell+1)}\tilde f_{k,\ell+1}
  \\&=d^h_{k+1,\ell}f'_{k+1,\ell}+d^v_{k,\ell+1}f'_{k,\ell+1}.
\end{align*}
We will use the above to
treat the cohomology of braided Hopf algebras with trivial
coefficients.

Let $\mathcal C$ be a braided monoidal category with braiding
$\braid$; we denote the unit object by $I$. We will assume that
the tensor product in $\mathcal C$ is strictly associative and
unital, so we can perform some calculations in the graphical
calculus (see for example \cite{Y} for more extensive examples of 
graphical calculations or \cite{JS} for more rigorous exposition). 
Our standard notations for braiding, multiplication, and unit will be
\begin{align*}
\braid=\braid_{XY}&=\gbeg23\got1X\got1Y\gnl\gbr\gnl\gob1Y\gob1X\gend,&
\mult=\mult_A&=\gbeg24\got1A\got1A\gnl\gcl1\gcl1\gnl\gmu\gnl\gob2A\gend,&
\eta=\eta_A&=\gbeg14\gnl\gu1\gnl\gcl1\gnl\gob1A\gend.
\end{align*}
Thus, if $A,B$ are algebras in $\mathcal C$, their tensor product
$A\ot B$ is defined to have multiplication
\[\mult_{A\ot B}=\gbeg44
                  \got1A\got1B\got1A\got1B\gnl
                  \gcl1\gbr\gcl1\gnl
                  \gmu\gmu\gnl
                  \gob2A\gob2B\gend.\]
We will need the following universal property of the tensor
product algebra: Given an algebra $R$ in $\mathcal C$ and algebra
morphisms $f_X\colon X\to R$ for $X\in\{A,B\}$ satisfying
$\mult_R(f_B\ot f_A)=\mult_R(f_A\ot f_B)\braid$, there exists a
unique algebra morphism $f\colon A\ot B\to R$ with
$f_A=f(1_A\ot\eta_B)$ and $f_B=f(\eta_A\ot 1_B)$, namely
$f=\mult_R(f_A\ot f_B)$.

We will denote by $A\olot B$ the tensor product algebra of
$A,B\in\mathcal C$ taken with respect to the inverse of the
braiding (i.e.\ in the braided category $(\mathcal
B,\braid\inv)$.). It possesses an obvious analogous universal
property. The equation
\[\gbeg46
  \got1B\got1A\got1B\got1A\gnl
  \gbr\gbr\gnl
  \gcl1\gbr\gcl1\gnl
  \gmu\gmu\gnl
  \gcn2122\gcn2122\gnl
  \gob2A\gob2B\gend
  =
  \gbeg46
  \got1B\got1A\got1B\got1A\gnl
  \gcl1\gibr\gcl1\gnl
  \gmu\gmu\gnl
  \gbbrh4226\gnl
  \gnl
  \gob2A\gob2B\gend
  \]
shows that $\braid\colon B\olot A\to A\ot B$ is an isomorphism of
algebras in $\mathcal C$.

The following lemma is a version of the Eckmann-Hilton argument
for algebras in braided monoidal categories.

\begin{lem}\label{two-mult}
Let $(A,\mult)$ be an algebra in the braided monoidal category
$\mathcal C$, equipped with a second multiplication $\multwo\colon
A\ot A\to A$ that shares the same unit $\eta\colon I\to A$.
Suppose that $\multwo\colon A\ot A\to A$ is multiplicative with
respect to $\mult$. Then the two multiplications coincide, and
they are commutative.
\end{lem}
\begin{proof}
We use the graphical symbols
\begin{align*}
\mult&=\gbeg24\got1A\got1A\gnl\gcl1\gcl1\gnl\gmu\gnl\gob2A\gend&
\multwo&=\gbeg24\got1A\got1A\gnl\gcl1\gcl1\gnl\gwmuc2\gnl\gob2A\gend&
\end{align*}
to distinguish the two multiplications.

The condition that $\multwo$ is multiplicative then reads
\[\gbeg46
  \got1A\got1A\got1A\got1A\gnl
  \gcl1\gcl1\gcl1\gcl1\gnl
  \gwmuc2\gwmuc2\gnl
  \gwmuh426\gnl
  \gcn4144\gnl
  \gob4A\gend
  =
  \gbeg46
  \got1A\got1A\got1A\got1A\gnl
  \gcl1\gbr\gcl1\gnl
  \gmu\gmu\gnl
  \gcn2122\gcn2122\gnl
  \gwmuch426\gnl
  \gob4A\gend\]
which we use twice in the calculation

\[\gbeg26
  \got 1A\got 1A\gnl
  \gbr\gnl
  \gcl1\gcl1\gnl
  \gwmuc2\gnl
  \gcn2122\gnl
  \gob2A\gend
  =
  \gbeg46
  \gvac1\got1A\got1A\gnl
  \gu1\gcl1\gcl1\gu1\gnl
  \gcl1\gbr\gcl1\gnl
  \gmu\gmu\gnl
  \gwmuch426\gnl
  \gob4A\gend
  =
  \gbeg46
  \gvac1\got1A\got1A\gnl
  \gu1\gcl1\gcl1\gu1\gnl
  \gwmuc2\gwmuc2\gnl
  \gwmuh426\gnl
  \gcn4144\gnl
  \gob4A\gend
  =
  \gbeg25
  \got1A\got1A\gnl
  \gcl1\gcl1\gnl
  \gmu\gnl
  \gcn2122\gnl
  \gob2A\gend
  =
  \gbeg46
  \got1A\gvac2\got1A\gnl
  \gcl1\gu1\gu1\gcl1\gnl
  \gwmuc2\gwmuc2\gnl
  \gwmuh426\gnl
  \gcn4144\gnl
  \gob4A\gend
  =
  \gbeg46
  \got1A\gvac2\got1A\gnl
  \gcl1\gu1\gu1\gcl1\gnl
  \gcl1\gbr\gcl1\gnl
  \gmu\gmu\gnl
  \gwmuch426\gnl
  \gob4A\gend
  =
  \gbeg25
  \got1A\got1A\gnl
  \gcl1\gcl1\gnl
  \gwmuc2\gnl
  \gcn2122\gnl
  \gob2A\gend\]
\end{proof}
\begin{rem}{\em
  It is not used in the proof that $(A,\mult)$ is an
  (associative) algebra. It would suffice to require that we are
  given two multiplications that share the same unit, and that one
  is multiplicative with respect to the other according to the
  same formula that results from the definition of a tensor
  product of algebras.}
\end{rem}

The above generalities on tensor products of algebras in a braided
category will be applied below to an $\N$-graded algebra that
occurs as the cohomology of an algebra. To fix the terminology, we
will denote the monoidal category of $\N$-graded objects in
$\mathcal C$ by $\mathcal C^\N$, with braiding $\braid_{\gr}$, and
we will call a commutative algebra in $\mathcal C^\N$ a braided
graded commutative algebra.

\begin{defn}{\em
  Let $\mathcal C$ be a monoidal category, and $A$ an
  augmented algebra in $\mathcal C$. The simplicial object $S_{\bullet}A$ is
  defined by $S_nA=A^{\ot n}$ with the faces
  $\partial_k^n\colon S_nA\to S_{n-1}A$ defined by $\partial _0^n=\varepsilon\ot (1_A)^{\ot (n-1)}$,
  $\partial _i^n=(1_A)^{\ot (i-1)}\ot\mult\ot (1_A)^{\ot (n-i-1)}$ for $0<i<n$ and
  $\partial _n^n=(1_A)^{\ot (n-1)}\ot\varepsilon$ and the degeneracies
  $\sigma_k^n\colon S_nA\to S_{n+1}A$ defined by inserting units in
  appropriate places.

  If $\mathcal C$ is abelian, then $S_{\bullet}A$ has an
  associated chain complex, which we will also denote by
  $S_\bullet A$.}
\end{defn}

Recall that a monoidal category $\mathcal C$ is (right) closed if
for each object $V$, the endofunctor $ - \ot V$ is left adjoint.
Its right adjoint is denoted $\ihom(V,\leer)$ and called an inner
hom-functor. If $\mathcal C$ is right closed, then $\ihom(V,Y)$ is
a bifunctor, covariant in $Y$ and contravariant in $V$.

More generally, one can denote by $\ihom(V,Y)$ an object, if it
exists, such that $\mathcal C(X\ot V,Y)\simeq\mathcal
C(X,\ihom(V,Y))$ as functors of $X\in\mathcal C$. By Yoneda's
Lemma, such hom-objects are unique, and bifunctorial in those
objects $V$ and $Y$ for which they exist.

Now consider a monoidal functor $\mathcal F\colon\mathcal
C\to\mathcal C'$. If $V,Y\in\mathcal C$ are such that the hom
objects $\ihom(V,Y)$ in $\mathcal C$ and $\ihom(\mathcal
FV,\mathcal FY)$ in $\mathcal C'$ exist, then there is a canonical
morphism $\mathcal F\ihom(V,Y)\to\ihom(\mathcal FV,\mathcal FY)$.
If this is an isomorphism, we say that $\mathcal F$ preserves the
hom object $\ihom(V,Y)$; if both categories are right closed, and
all hom objects are preserved, we say that $\mathcal F$ preserves
inner hom-functors.

\begin{defn}\label{HHdef} {\em
  Let $\mathcal C$ be an abelian monoidal category,
  and $A$ an augmented algebra in $\mathcal C$ such that all
  hom objects $\ihom(A^{\ot n},I)$ exist (e.g.\ if $\mathcal C$ is closed).

  The Hochschild cohomology of $A$ in
  $\mathcal C$ (with trivial coefficients)
  is the cohomology of the cochain complex
  $\ihom(S_{\bullet}A,I)$. Thus the Hochschild cohomology consists of objects in
  the category $\mathcal C$. We denote it by ${\HH}^*(A)$.}
\end{defn}

\begin{rem}
 {\em We will only need the results in this section for the case $\mathcal
 C=\YD$, but we will indicate more generally how they apply to
 ordinary algebras and (braided) Hopf algebras. For this, it is
 not necessary always to stick to the paper's general assumption
 that we work over a base \emph{field}, therefore we assume now
 that $k$ is an arbitrary commutative base ring.
 \begin{enumerate}
 \item  In the case of an ordinary
  augmented $k$-algebra $A$, Definition \ref{HHdef} recovers the ordinary definition of Hochschild
  cohomology with trivial coefficients, that is ${\HH}^*(A)={\HH}^*(A,k)$ (see \S 2), which
  motivates our notation; general coefficients are
  not considered here.
 \item Assume that $\mathcal F\colon\mathcal
  C\to\mathcal C'$ is an exact monoidal functor that preserves the
  hom objects $\ihom(S_\bullet A,I)$. Then $\mathcal F$
  preserves Hochschild cohomology of $A$ in the sense that
  ${\HH}^*(\mathcal F(A))\simeq \mathcal F({\HH}^*(A))$.
 \item The category of (say, left) $G$-modules for a
  $k$-Hopf algebra $G$ has inner hom-objects preserved by the underlying
  functor to the category of $k$-modules. Also, if $G$ is a finitely
  generated projective Hopf algebra over $k$, the categories of (say,
  left) $G$-comodules, and of
  Yetter-Drinfeld modules over $G$ have inner hom-objects that are
  preserved by the underlying functors to
  the category of $k$-modules. More concretely, the hom-object
  $\ihom(V,Y)$ is $\Hom_k(V,Y)$ with $G$-module and comodule
  structures induced by the antipode and its inverse.
 \item Recall that a dual object $(V^*,e,d)$ of an object $V$ is an
  object $V^*$ with morphisms $e\colon V^*\ot V\to I$ and $d\colon
  I\to V\ot V^*$ that satisfy $(1_V\ot e)(d\ot 1_V)=1_V$ and $(e\ot
  1_{V^*})(1_{V^*}\ot d)=1_{V^*}$.  If an object $V$ in $\mathcal C$
  has a dual, then all the inner hom-objects $\ihom(V,Y)$ exist and
  are given by $\ihom(V,Y)=Y\ot V^*$. Also, monoidal functors
  always preserve dual objects (by and large because these are
  given by morphisms and relations) and thus, if $V$ has a dual,
  then the inner hom-objects $\ihom(V,Y)$ are preserved by every
  monoidal functor.
 \item A module $V$ over a commutative ring $k$ has a dual in the
  category of $k$-modules if and only if
  it is finitely generated projective; then $V^*$ is the dual
  module, $e$ is the evaluation, and $d$ maps $1_k$ to the canonical
  element, in the field case obtained by tensoring basis elements
  with the elements of the dual basis. If $G$ is a $k$-Hopf algebra
  with bijective antipode, then in the categories of $G$-modules,
  $G$-comodules, and Yetter-Drinfeld modules over $G$ an object $V$
  has a dual if and only if $V$ is finitely generated projective
  over $k$.
 \item Let $G$ be a flat Hopf algebra over $k$ with bijective antipode, so
  that the category $\mathcal C$ of Yetter-Drinfeld modules over $G$ is abelian
  braided monoidal with an exact underlying functor to the category
  $\mathcal C'$ of $k$-modules. Then by the above remarks
  the underlying functor will preserve the Hochschild
  cohomology of $A$ in $\mathcal C$ whenever $G$ is finitely generated projective over
  $k$, or $A$ is finitely generated
  projective. Thus, in either of these cases the Hochschild cohomology of
  $A$ defined as above within the Yetter-Drinfeld category is the
  same as the ordinary Hochschild cohomology with trivial
  coefficients, endowed with a Yetter-Drinfeld module structure
  induced by that of $A$.
\end{enumerate}}
\end{rem}

\begin{rem}
\label{xi}{\em
  In a closed monoidal category $\mathcal C$, there is a natural
  morphism $\xi\colon\ihom(B,I)\ot\ihom(A,I)\to\ihom(A\ot B,I)$ which is
  suitably coherent with respect to higher tensor products. Under the adjunction
  \[\mathcal C(\ihom(B,I)\ot\ihom(A,I),\ihom(A\ot B,I))\simeq
  \mathcal C(\ihom(B,I)\ot\ihom(A,I)\ot A\ot B,I)\]
  $\xi$ corresponds to
  \[\ihom(B,I)\ot\ihom(A,I)\ot A\ot B\xrightarrow{1\ot e\ot 1}\ihom(B,I)\ot B\xrightarrow{e}I,\]
  where $e\colon\ihom(X,I)\ot X\to I$ is the natural
  ``evaluation'' morphism.
  Of course, to define the morphism in
  question, it is enough to assume that all the inner hom objects
  that occur in it exist.

  If all inner
  hom functors $\ihom(V,I)$ exist, the isomorphisms $\xi$ make
  $\ihom(\leer,I)$ into a contravariant weak monoidal
  functor.

  If $A$ and $B$ have dual objects in $\mathcal C$, then one can
  show that $\xi$ above is an isomorphism; this can be attributed
  to the fact that $(B^*\ot A^*,\tilde e,\tilde d)$ is a dual
  object for $A\ot B$ with
  \[\tilde e=\left(B^*\ot A^*\ot A\ot B\xrightarrow{1\ot e\ot 1}B^*\ot B\xrightarrow e I\right)\]
  \[\tilde d=\left(I\xrightarrow dA\ot A^*\xrightarrow{1\ot d\ot 1}A\ot B\ot B^*\ot
  A^*\right)\]}
\end{rem}

\begin{defn}\label{cupdef}{\em
  The cup product on ${\HH}^*(A)$ is the collection of morphisms
  ${\HH}^n(A)\ot{\HH}^m(A)\to {\HH}^{m+n}(A)$ induced by the morphisms
  \[\ihom(S_nA,I)\ot\ihom(S_mA,I)\xrightarrow\xi\ihom(S_mA\ot S_nA,I)= \ihom(S_{m+n}(A),I).\]}
\end{defn}
\begin{rem}{\em
Thus, in the category of modules over a commutative base ring $k$,
our definition of cup product recovers \emph{the opposite} of the
usual cup product in Hochschild cohomology.
}
\end{rem}

{F}rom now on we assume that $\mathcal C$ is braided.

\begin{lem}
  Let $A,B$ be two augmented algebras in $\mathcal C$. Denote by
  $S_\bullet A\times S_\bullet B$ the bisimplicial object obtained
  by tensoring the two simplicial objects $S_\bullet A$ and
  $S_\bullet B$. An isomorphism
  $S_\bullet(A\ot B)\to\diag(S_\bullet A\times S_\bullet B)$ of
  simplicial objects is given by the morphisms
  \[g_n\colon S_n(A\ot B)=(A\ot B)^{\ot n}\to A^{\ot n}\ot B^{\ot
  n}=\diag_n(S_\bullet A\ot S_\bullet B)\]
  that are composed from instances of the braiding
  $\braid_{B,A}\colon B\ot A\to A\ot B$.
\end{lem}
\begin{proof}
  For example we have
  \[g_3=\gbeg64
        \got1A\got1B\got1A\got1B\got1A\got1B\gnl
        \gcl2\gbr\gbr\gcl2\gnl
        \gvac1\gcl1\gbr\gcl1\gnl
        \gob1A\gob1A\gob1A\gob1B\gob1B\gob1B\gend
        \quad\text{and}\quad
    g_4=\gbeg85
        \got1A\got1B\got1A\got1B\got1A\got1B\got1A\got1B\gnl
        \gcl3\gbr\gbr\gbr\gcl3\gnl
        \gvac1\gcl2\gbr\gbr\gcl2\gnl
        \gvac2\gcl1\gbr\gcl1\gnl
        \gob1A\gob1A\gob1A\gob1A\gob1B\gob1B\gob1B\gob1B\gend\]
  and therefore
  \[(\partial_2\ot\partial_2)g_4
    =\gbeg87
     \got1A\got1B\got1A\got1B\got1A\got1B\got1A\got1B\gnl
     \gcl5\gbr\gbr\gbr\gcl5\gnl
     \gvac1\gcl2\gbr\gbr\gcl2\gnl
     \gvac2\gcl1\gbr\gcl1\gnl
     \gvac1\gmu\gcl2\gcl2\gmu\gnl
     \gvac1\gcn2122\gvac2\gcn2122\gnl
     \gob1A\gob2A\gob1A\gob1B\gob2B\gob1B\gend
     =
     \gbeg88
     \got1A\got1B\got1A\got1B\got1A\got1B\got1A\got1B\gnl
     \gcl6\gcl2\gcl1\gbr\gcl1\gcl2\gcl6\gnl
     \gvac2\gmu\gmu\gnl
     \gvac1\gbbrh3214\gbbrh3225\gnl
     \gnl
     \gvac1\gcl2\gbbrh4226\gcl2\gnl
     \gnl
     \gob1A\gob1A\gob2A\gob2B\gob1B\gob1B\gend
     =g_3\partial_2
  \]
  The general claim $(\partial_k\ot\partial_k)g_n=g_{n-1}\partial_k$
  is a merely larger version of the above
  example for $0<k<n$, and simpler (with augmentations instead of
  multiplications) for $k\in\{0,n\}$.
\end{proof}

According to the discussion at the beginning of the section we
have two chain homotopic morphisms
\[\diag(S_\bullet A\times S_\bullet B)\simeq S_\bullet A\ot S_\bullet B.\]
Here, the right hand side denotes the tensor product of chain
complexes. Explicitly, the two morphisms are given
by
\begin{gather*}
  f_{k\ell}=(1_A)^{\ot k}\ot\varepsilon^{\ot\ell}\ot \varepsilon^{\ot k}\ot (1_B)^{\ot\ell}\colon S_nA\ot S_nB\to S_kA\ot S_\ell B
  \\
  f'_{k\ell}=(-1)^{k\ell}\varepsilon^{\ot\ell}\ot (1_A)^{\ot k}\ot (1_B)^{\ot \ell}\ot\varepsilon^{\ot
  k}\colon S_nA\ot S_nB\to S_kA\ot S_\ell B
\end{gather*}
for $k+\ell=n$.

Upon application of $\ihom(\leer,I)$, each of these morphisms,
composed with the map $g^{-1}$ from the preceding lemma, yields a
morphism
\begin{equation}\label{Tprime}
\ihom(S_\bullet B,I)\ot\ihom(S_\bullet A,I)\to\ihom(S_\bullet
A\ot S_\bullet B,I)\to\ihom(S_\bullet(A\ot B),I).
\end{equation} Both these maps, being homotopic, yield the same morphism
\[T'\colon{\HH}^\bullet(B)\ot{\HH}^\bullet(A)\to{\HH}^\bullet(A\ot B),\]
which we compose with the inverse of the braiding in $\mathcal C$
to obtain a morphism
\[T=T'\braid_{\gr}\inv\colon{\HH}^\bullet(A)\ot{\HH}^\bullet(B)\to{\HH}^\bullet(A\ot B).\]

\begin{lem}\label{morphism-T}
  The morphism
  \[T\colon{\HH}^\bullet(A)\ot {\HH}^\bullet(B)\to{\HH}^\bullet(A\ot B)\]
  is a morphism of graded algebras in $\mathcal C$.
\end{lem}
\begin{proof}
  We need to show that
  \[T'\colon{\HH}^\bullet(B)\olot{\HH}^\bullet(A)\to{\HH}^\bullet(A\ot B)\]
  is a morphism of algebras. To do this, we
  consider the algebra morphisms
  \[{\HH}^\bullet(A)\to{\HH}^\bullet(A\ot B)\leftarrow{\HH}^\bullet(B).\]
  We need to show that the composite
  \[{\HH}^\bullet(A)\ot{\HH}^\bullet(B)\to{\HH}^\bullet(A\ot B)\ot{\HH}^\bullet(A\ot B)\to{\HH}^\bullet(A\ot B),\]
  where the last morphism is multiplication, is the map in the
  statement of the lemma, while the composite
  \[{\HH}^\bullet(B)\ot{\HH}^\bullet(A)\to{\HH}^\bullet(A\ot B)\ot{\HH}^\bullet(A\ot B)\to{\HH}^\bullet(A\ot B)\]
  is the same, composed with the inverse braiding. But for this we
  only need to observe that
  \[S_n(A\ot B)=S_k(A\ot B)\ot S_\ell(A\ot B)\xrightarrow{S_k(1_A\ot\varepsilon)\ot S_\ell(\varepsilon\ot 1_B)} S_k(A)\ot S_\ell(B)\]
  equals
  \[S_n(A\ot B)\xrightarrow{g} S_n(A)\ot S_n(B)\xrightarrow{f_{k\ell}} S_k(A)\ot S_\ell(B),\]
  whereas \[S_n(A\ot B)=S_k(A\ot B)\ot S_\ell(A\ot
  B)\xrightarrow{S_k(\varepsilon\ot 1_B)\ot S_\ell(1_A\ot\varepsilon)} S_k(B)\ot S_\ell(A)\]
  equals
  \[S_n(A\ot B)\xrightarrow{g} S_n(A)\ot S_n(B)\xrightarrow{f'_{\ell k}}
  S_\ell(A)\ot S_k(B)\xrightarrow{\braid_{\gr}^{-1}} S_k(B)\ot S_\ell(A).\]
\end{proof}

\begin{remark}{\em
In the case which is of main interest to us,  $k$ is a field, $\mathcal C$ 
is the category of Yetter-Drinfeld modules over a Hopf algebra $G$, 
and $A, B$ are finite dimensional. We briefly sketch why these conditions imply 
that $T$ is an isomorphism (see also \cite[Prop.4.5]{GM}). 

First, both in the construction of $T'$ (when homology is applied to the sequence \ref{Tprime}) 
and in that of the cup product in Definition \ref{cupdef} we have tacitly used a 
version of a K\"unneth-type morphism 
\[\HH(C^\bullet)\ot\HH(D^\bullet)\to \HH(C^\bullet\ot D^\bullet)\]
for chain complexes $C^\bullet,D^\bullet$ in a closed monoidal category $\mathcal C$. 
Under suitable flatness hypotheses on the objects and homologies that
occur, this map will be an isomorphism. For example, this is
the case if the tensor product in $\mathcal C$ is assumed to be
exact in each argument which holds in our situation.

Second, the morphism $\xi$ that occurs in the construction of $T'$ will be 
an isomorphism by Remark \ref{xi} since the algebras $A$ and $B$ clearly have 
duals in $\mathcal C$. 

Therefore, in our case $T$ is an isomorphism.
}\end{remark}

\begin{thm}
Let $R$ be a bialgebra in the abelian braided monoidal category
$\mathcal C$. Then Hochschild cohomology ${\HH}^*(R)$ is a braided
graded commutative algebra in $\mathcal C$, assuming it is defined
(i.e.~the necessary hom-objects exist).
\end{thm}
\begin{proof}
Comultiplication $\Delta$ on the bialgebra $R$ yields a second
multiplication
\[\multwo:=\left({\HH}^*(R)\ot{\HH}^*(R)\to{\HH}^*(R\ot R)\xrightarrow{{\HH}^*(\Delta)}{\HH}^*(R)\right)\]
on ${\HH}^*(R)$. By construction, $\multwo$ is an algebra map with
respect to the algebra structure of ${\HH}^*(R)$ given by the cup
product. Also, the two multiplications share the same unit,
induced by the counit of $R$. Therefore, the two multiplications
coincide and are commutative by Lemma \ref{two-mult}.
\end{proof}

\begin{cor}\label{cor:braided-graded}
Let $G$ be a Hopf algebra with bijective antipode, and $R$ a
bialgebra in the category $\mathcal C$ of Yetter-Drinfeld modules
over $G$. Assume that either $G$ or $R$ is finite dimensional.
Then the Hochschild cohomology of $R$ in $\mathcal C$ is a braided
graded commutative algebra. Consequently, the opposite of ordinary
Hochschild cohomology is a braided commutative Yetter-Drinfeld
module algebra.
\end{cor}
\begin{proof}
In fact the assumptions on $G$ or $R$ are designed to ensure that
the hypotheses on inner hom objects in the results above are
satisfied: If $G$ is finitely generated projective, the category
of Yetter-Drinfeld modules as a whole is closed, with inner
hom-functors preserved by the underlying functors. If $G$ is not
finitely generated projective, but $R$ is, then still $R$ and all
its tensor powers have dual objects in the category of
Yetter-Drinfeld modules over $G$, that are preserved by the
underlying functors.
\end{proof}

\section{Cohomology of quantum complete intersections}\label{sec-step0}

Let $\theta$ be a positive integer, and for each $i$, $1\leq i\leq \theta$,
let $N_i$ be an integer greater than 1.
Let $q_{ij}\in k^{\times}$ for $1\leq i<j\leq \theta$.
Let $S$ be the $k$-algebra generated by $x_1,\ldots,x_{\theta}$ subject
to the relations
\begin{equation}\label{qci}
   x_ix_j=q_{ij}x_jx_i \mbox{ for all } i<j \  \ \
   \mbox{ and } \ \ \ x_i^{N_i}=0 \mbox{ for all } i.
\end{equation}
It is convenient to set $q_{ji}=q_{ij}^{-1}$ for $i<j$.
For ${\mathcal B}(V)$, the Nichols algebra constructed
in Section \ref{pointed},
Lemma \ref{filtration} implies that $\Gr {\mathcal B}(V)$ has this form.

We will compute $\coh^*(S,k) =\Ext^*_S(k,k)$ for use in later
sections.
The structure of this ring may be determined by using the
braided K\"unneth formula
of Grunenfelder and Mastnak \cite[Prop.\ 4.5 and Cor.\ 4.7]{GM}
or the twisted tensor product formula of Bergh and Oppermann
\cite[Thms.\ 3.7 and 5.3]{BO}.
We give details, using an explicit free $S$-resolution of $k$, in order
to record needed information at the chain level.
Our resolution was originally adapted from \cite{BG} (see also \cite{HK}):
It is a braided tensor product of the periodic resolutions
\begin{equation}\label{perres}
\cdots \!\stackrel{x_i^{N_i-1}\cdot}{\relbar\joinrel
\longrightarrow}\! k[x_i]/(x_i^{N_i})\!
\stackrel{x_i\cdot}{
\longrightarrow}\! k[x_i]/(x_i^{N_i})\!
\stackrel{x_i^{N_i-1}\cdot}{\relbar\joinrel
\longrightarrow} \! k[x_i]/(x_i^{N_i})\!
\stackrel{x_i\cdot}{
\longrightarrow} \! k[x_i]/(x_i^{N_i})\! \stackrel{\varepsilon}
{\rightarrow}\! k\rightarrow \! 0,
\end{equation}
one for each $i$, $1\leq i\leq \theta$.

Specifically, let $K_{\bullet}$ be the following complex of free
$S$-modules. For each $\theta$-tuple $(a_1,\ldots,a_{\theta})$ of
nonnegative integers, let $\Phi(a_1,\ldots,a_{\theta})$ be a free
generator in degree $a_1+\cdots +a_{\theta}$. Thus
$K_n=\oplus_{a_1+\cdots +a_{\theta}=n}
S\Phi(a_1,\ldots,a_{\theta})$. For each $i$, $1\leq i\leq\theta$,
let $\sigma_i, \tau_i: \N\rightarrow \N$ be the functions defined
by
$$
  \sigma_i(a) = \left\{\begin{array}{cl}
              1, & \mbox{ if }a\mbox{ is odd}\\
            N_i-1, & \mbox{ if }a\mbox{ is even},
              \end{array}\right.
$$
and $\tau_i(a) = \displaystyle{\sum_{j=1}^a \sigma_i(j)}$ for $a\geq 1$, $\tau(0)=0$.
Let
$$
  d_i (\Phi(a_1,\ldots,a_{\theta})) = \left(\prod_{\ell < i} (-1)^{a_{\ell}}
     q_{\ell i} ^{\sigma_i(a_i)\tau_{\ell}(a_{\ell})}\right)
       x_i^{\sigma_i(a_i)} \Phi(a_1,\ldots, a_i -1,\ldots , a_{\theta})
$$
if $a_i>0$, and $d_i (\Phi(a_1,\ldots,a_{\theta})) =0$ if $a_i=0$.
Extend each $d_i$ to an $S$-module homomorphism.
Note that $d_i^2=0$ for each $i$ since $x_i^{N_i}=0$ and
$\sigma_i(a_i)+\sigma_i(a_i -1)=N_i$.
If $i<j$, we have

\smallskip

\noindent
$
  d_id_j(\Phi(a_1,\ldots,a_{\theta}))
$
\begin{eqnarray*}
  &=&        d_i \left(\left(\prod_{\ell <j} (-1)^{a_{\ell}} q_{\ell j}^{\sigma_j(a_j)
     \tau_{\ell}(a_{\ell})}\right) x_j^{\sigma_j(a_j)} \Phi(a_1,\ldots,
        a_j -1, \ldots,a_{\theta})\right)\\
   &=& \left(\prod_{\ell <j} (-1)^{a_{\ell}} q_{\ell j}^{\sigma_j(a_j)
  \tau_{\ell}(a_{\ell})}\right) \left(\prod_{m<i} (-1)^{a_m} q_{mi}^{\sigma_i(a_i)
    \tau_m(a_m)}\right)\cdot\\
    &&\hspace{1in} x_j^{\sigma_j(a_j)} x_i^{\sigma_i(a_i)}
   \Phi(a_1,\ldots, a_i-1,\ldots,a_j-1,\ldots,a_{\theta}).
\end{eqnarray*}
Now
$d_jd_i(\Phi(a_1,\ldots,a_{\theta}))$ is also an $S$-multiple of
$\Phi(a_1,\ldots,a_i-1,\ldots, a_j-1,\ldots,a_{\theta})$;
comparing the two, we see that the term in which $\ell =i$ has a
scalar factor that changes from
$(-1)^{a_i} q_{ij}^{\sigma_j(a_j)\tau_i(a_i)}$ to
$(-1)^{a_i-1} q_{ij}^{\sigma_j(a_j) \tau_i(a_i -1)}$, and
$x_j^{\sigma_j(a_j)} x_i^{\sigma_i(a_i)}$ is replaced by
$x_i^{\sigma_i(a_i)}x_j^{\sigma_j(a_j)}=q_{ij}^{\sigma_j(a_j)
    \sigma_i(a_i)} x_j^{\sigma_j(a_j)}x_i^{\sigma_i(a_i)}$.
Since $\tau_i(a_i)=\tau_i(a_i -1)+\sigma_i(a_i)$, this shows that
$$
  d_id_j + d_jd_i =0.
$$
Letting
\begin{equation}\label{dsum}
d=d_1+\cdots +d_{\theta},
\end{equation}
we now have
$d^2=0$, so $K_{\bullet}$ is indeed a complex.

Next we show that $K_{\bullet}$ is a resolution of $k$ by giving a contracting
homotopy:
Let $\alpha\in S$, and fix $\ell$, $1\leq \ell\leq \theta$. Write
$\alpha =\sum_{j=0}^{N_i-1}\alpha_jx_{\ell}^j$
where $\alpha_j$ is in the subalgebra of $S$ generated by the $x_i$ with
$i\neq \ell$. Define $s_{\ell}(\alpha\Phi(a_1,\ldots, a_{\theta}))$ to be
the sum $\sum_{j=0}^{N_i -1} s_{\ell}(\alpha_j x_{\ell}^j \Phi(a_1,\ldots, a_{\theta}))$,
where

\smallskip

\noindent
$
   s_{\ell}(\alpha_jx_{\ell}^j\Phi(a_1,\ldots,a_{\theta}))
$
$$=
    \left\{ \begin{array}{l}
   \delta_{j>0} \displaystyle{\left(\prod_{m<\ell}(-1)^{a_m} q_{m\ell}^{-\sigma_{\ell}(a_{\ell}+1)
   \tau_m(a_m)} \right)}
     \alpha_j x_{\ell}^{j-1}\Phi(a_1,\ldots,a_{\ell}+1,\ldots,a_{\theta}),\\
    \hspace{3in} \mbox{ if } a_{\ell} \mbox{ is even}\\
    \delta_{j,N_{\ell}-1}\displaystyle{\left( \prod_{m<\ell} (-1)^{a_m} q_{m\ell}^{-\sigma_{\ell}(a_{\ell}+1)
    \tau_m(a_m)}\right)} \alpha_j\Phi(a_1,\ldots,a_{\ell}+1,\ldots,a_{\theta}),\\
     \hspace{3in} \mbox{ if } a_{\ell} \mbox{ is odd}
   \end{array}\right.
$$
where $\delta_{j>0} =1$ if $j>0$ and $0$ if $j=0$.
Calculations show that for all $i$, $1\leq i\leq\theta$,
$$
   (s_id_i + d_is_i)(\alpha_jx_i^j\Phi(a_1,\ldots,a_{\theta}))  =
   \left\{\begin{array}{cl}\alpha_jx_i^j\Phi(a_1,\ldots,a_{\theta}), & \mbox{if }j>0\mbox{ or }a_i>0\\
        0, & \mbox{if }j=0\mbox{ and } a_i=0\end{array}\right.
$$
and $s_{\ell}d_i+d_is_{\ell} = 0$  for all $i,\ell$ with $i\neq \ell$.
For each $x_1^{j_1}\cdots x_{\theta}^{j_{\theta}}\Phi(a_1,\ldots,a_{\theta})$,
let $c = c_{j_1,\ldots,j_{\theta},a_1,\ldots,a_{\theta}}$ be the cardinality
of the set of all $i$ ($1\leq i\leq \theta$)
such that $j_ia_i=0$.
Define
$$
   s(x_1^{j_1}\cdots x_{\theta}^{j_{\theta}}\Phi(a_1,\ldots,a_{\theta}))
   = \frac{1}{\theta - c_{j_1,\ldots,j_{\theta},a_1,\ldots,a_{\theta}}}
   (s_1+\cdots +s_{\theta})
   (x_1^{j_1}\cdots x_{\theta}^{j_{\theta}}\Phi(a_1,\ldots,a_{\theta})).
$$
Then $sd+ds=\id$ on each $K_n$, $n>0$.
That is, $K_{\bullet}$ is exact in positive degrees.
To show that $K_{\bullet}$ is a resolution of $k$, put $k$ in degree $-1$,
and let the corresponding differential be the
counit map $\varepsilon: S\rightarrow k$.
It can be shown directly that $K_{\bullet}$ then becomes exact at $K_0=S$:
The kernel of $\varepsilon$ is spanned over the field $k$ by the elements
$x_1^{j_1}\cdots x_{\theta}^{j_{\theta}} \Phi(0,\ldots,0)$, $0\leq
j_i\leq N_i$, with at least one $j_i\neq 0$.
Let $x_1^{j_1}\cdots x_{\theta}^{j_{\theta}} \Phi(0,\ldots,0)$ be such
an element, and let $i$ be the smallest positive integer such that $j_i\neq
0$.
Then
  $d(x_i^{j_i-1}\cdots x_{\theta}^{j_{\theta}}\Phi(0,\ldots, 1,\ldots, 0))$
is a nonzero scalar multiple of
  $x_i^{j_i}\cdots x_{\theta}^{j_{\theta}} \Phi(0,\ldots,0)$.
Thus $\ker (\varepsilon) = \im (d)$, and $K_{\bullet}$ is a free resolution
of $k$ as an $S$-module.

Next we will use $K_{\bullet}$ to compute $\Ext^*_S(k,k)$.
Applying $\Hom_S( - , k)$ to $K_{\bullet}$, the induced differential $d^*$
is the zero map since $x_i^{\sigma_i(a_i)}$ is always in the
augmentation ideal. Thus the cohomology {\em is} the complex
$\Hom_S(K_{\bullet},k)$, and in degree $n$ this is a vector space of
dimension $\binom{n+\theta -1}{\theta -1}$.
Now let $\xi_i\in\Hom_S(K_2,k)$ be the function dual to
$\Phi(0,\ldots,0,2,0,\ldots, 0)$ (the 2 in the $i$th place)
and $\eta_i\in \Hom_S(K_1,k)$ be the function dual to
$\Phi(0,\ldots,0,1,0,\ldots,0)$ (the 1 in the $i$th place).
By abuse of notation, identify these functions with the
corresponding elements in $\coh^2(S,k)$ and $\coh^1(S,k)$, respectively.
We will show that the $\xi_i$, $\eta_i$ generate $\coh^*(S,k)$, and
determine the relations among them. In order to do this we will
abuse notation further and denote by $\xi_i$ and $\eta_i$ the
corresponding chain maps $\xi_i : K_n\rightarrow K_{n-2}$
and $\eta_i : K_n \rightarrow K_{n-1}$ defined by
\begin{eqnarray*}
  \xi_i(\Phi(a_1,\ldots,a_{\theta}))\!\!\! &=& \!\!\! \prod_{\ell>i} q_{i\ell}
    ^{N_i\tau_{\ell}(a_{\ell})}\Phi(a_1,\ldots,a_i-2,\ldots,a_{\theta})\\
  \eta_i (\Phi(a_1,\ldots,a_{\theta}))
 \!\!\! &=& \!\!\! \prod_{\ell<i}q_{\ell i} ^{(\sigma_i(a_i)-1)\tau_{\ell}
   (a_{\ell})} \prod_{\ell >i} (-1)^{a_{\ell}} q_{i\ell}^{\tau_{\ell}
  (a_{\ell})} x_i^{\sigma_i(a_i)-1} \Phi(a_1,\ldots,a_i-1,\ldots,a_{\theta}).
\end{eqnarray*}
Calculations show that these are indeed chain maps.
The ring structure of the subalgebra of $\coh^*(S,k)$
generated by $\xi_i,\eta_i$ is given by composition of these chain
maps. Direct calculation shows that the relations
given in Theorem \ref{step0} below hold.
(Note that if $N_i\neq 2$ the last relation implies
$\eta_i^2=0$.)
Alternatively, in case $S={\mathcal B}(V)$, a Nichols algebra
defined in Section \ref{pointed}, we may apply Corollary \ref{cor:braided-graded}
to obtain the relations in Theorem \ref{step0} below,
performing calculations that are a special case of the ones in
the proof of Theorem \ref{relns}.
Thus any element in the algebra generated by the $\xi_i$ and $\eta_i$
may be written as a linear combination of elements of the form
$\xi_1^{b_1}\cdots \xi_{\theta}^{b_{\theta}}\eta_1^{c_1}\cdots
\eta_{\theta}^{c_{\theta}}$ with $b_i\geq 0$ and $c_i\in\{0,1\}$.
Such an element takes $\Phi(2b_1+c_1,\ldots,2b_{\theta}+c_{\theta})$ to
a nonzero scalar multiple of $\Phi(0,\ldots,0)$ and all other
$S$-basis elements of $K_{\sum (2b_i+c_i)}$ to 0.
Since the dimension of $\coh^n(S,k)$ is $\binom{n+\theta -1}{\theta -1}$,
this shows that
the $\xi_1^{b_1}\cdots \xi_{\theta}^{b_{\theta}}\eta_1^{c_1}\cdots
\eta_{\theta}^{c_{\theta}}$ form a $k$-basis for $\coh^*(S,k)$.
We have proven:

\begin{thm}\label{step0}
Let $S$ be the $k$-algebra generated by $x_1,\ldots,x_{\theta}$,
subject to relations (\ref{qci}). Then
$\coh^*(S,k)$ is generated by $\xi_i, \ \eta_i$ ($i=1,\ldots,\theta$)
where $\deg \xi_i = 2$ and $\deg \eta_i = 1$,
subject to the relations
\begin{equation}\label{relations}
  \xi_i\xi_j=q_{ji}^{N_iN_j}\xi_j\xi_i, \ \ \
  \eta_i\xi_j=q_{ji}^{N_j}\xi_j\eta_i, \ \ \
   \mbox{ and } \ \eta_i\eta_j=-q_{ji}\eta_j\eta_i.
\end{equation}
\end{thm}

Note that although the relations (\ref{relations}) can be obtained as
a consequence of Corollary \ref{cor:braided-graded}, in order to obtain the full
statement of the theorem, we needed more information.

\begin{remark}{\em
We obtain \cite[Prop.\ 2.3.1]{GK} as a corollary: In this case we
replace the generators $x_i$ of $S$ by the generators
$E_{\alpha}$ ($\alpha\in\Delta^+$) of $\Gr u_{q}^+$, whose
relations are
$$
  E_{\alpha}E_{\beta}=q^{\langle\alpha,\beta\rangle} E_{\beta}
  E_{\alpha} \ (\alpha\succ \beta) , \ \ \ \ \ \
  E_{\alpha}^{\ell} = 0 \ (\alpha\in\Delta^+).
$$
Theorem \ref{step0} then implies that $\coh^*(\Gr u_{q}^+)$ is generated
by $\xi_{\alpha},\eta_{\alpha}$ ($\alpha\in\Delta^+$) with relations
$$
  \xi_{\alpha}\xi_{\beta}=\xi_{\beta}\xi_{\alpha}, \ \
   \eta_{\alpha}\xi_{\beta}=\xi_{\beta}\eta_{\alpha}, \ \
  \eta_{\alpha}\eta_{\beta}=-q^{-\langle\alpha,\beta\rangle}
  \eta_{\beta}\eta_{\alpha} \ (\alpha\succ\beta), \ \
  \eta_{\alpha}^2=0,
$$
which is precisely \cite[Prop.\ 2.3.1]{GK}.}\end{remark}

Now assume a finite group $\Gamma$ acts on $S$ by automorphisms
for which $x_1,\ldots, x_{\theta}$ are eigenvectors.
For each $i$, $1\leq i\leq \theta$, let $\chi_i$ be the character on $\Gamma$
for which $g x_i = \chi_i(g)x_i$.
As the characteristic of $k$ is 0, we have
$$
  \Ext^*_{S\# k\Gamma}(k,k) \simeq \Ext^*_S(k,k)^{\Gamma},
$$
where the action of $\Gamma$ at the chain level on (\ref{perres})
is as usual in
degree 0, but shifted in higher degrees so as to make the differentials
commute with the action of $\Gamma$.
Specifically,
note that the following action of $\Gamma$ on $K_{\bullet}$
commutes with the differentials:
$$
  g\cdot \Phi(a_1,\ldots,a_{\theta})=\prod_{\ell =1}^{\theta}
   \chi_{\ell}(g)^{\tau_{\ell}(a_{\ell})}\Phi(a_1,\ldots,a_{\theta})
$$
for all $g\in \Gamma$, and $a_1,\ldots,a_{\theta}\geq 0$.
Then  the induced action of $\Gamma$  on generators $\xi_i, \eta_i$ of the
cohomology ring $\coh^*(S,k)$ is given explicitly by
\begin{equation}
\label{action}
g\cdot\xi_i=\chi_i(g)^{-N_i}\xi_i \text{ and }  g\cdot\eta_i =\chi_i(g)^{-1}\eta_i.
\end{equation}


\section{Coradically graded finite dimensional pointed Hopf algebras}\label{graded}
Let ${\mathcal D}$ be arbitrary data as in (\ref{data}),
$V$ the corresponding Yetter-Drinfeld module,
and $R={\mathcal B}(V)$ its Nichols algebra, and described in Section \ref{pointed}.
By Lemma \ref{filtration}, there is a filtration on $R$ for which
$S=\Gr R$ is of type $A_1\times \cdots\times A_1$, given by generators and
relations of type (\ref{qci}).
Thus $\coh^*(S,k)$ is given by
Theorem \ref{step0}.  As the filtration is finite, there is  a convergent spectral
sequence associated to the filtration
(see \cite[5.4.1]{W}):
\begin{equation}
\label{sseq1}
\xymatrix
{E^{p,q}_1 = \coh^{p+q}(\Gr_{(p)} R, k) \ar@{=>}[r] & \coh^{p+q}(R,k) .
}
\end{equation}
It follows that the $E_1$-page  of the spectral sequence  is given by
Theorem~\ref{step0} with grading corresponding to the filtration on $R$.
We will see that by (\ref{xialpha}) and Lemma \ref{perm} below,
the generators $\xi_i$ are in degrees  $(p_i, 2-p_i)$,
where
\begin{equation}\label{pi}
  p_i = N_{\beta_1}\cdots N_{\beta_i}(N_{\beta_i}\cdots
  N_{\beta_r}\Ht(\beta_i)+1).
\end{equation}
Since the PBW basis elements (\ref{PBW}) are eigenvectors for $\Gamma$,
the action of $\Gamma$ on $R$ preserves the filtration, and
we further get a spectral sequence  converging to the cohomology of
$u({\mathcal D},0,0)\simeq R\# k\Gamma$:
\begin{equation}
\label{sseq2}
\xymatrix
{\coh^{p+q}(\Gr_{(p)} R, k)^\Gamma \ar@{=>}[r] & \coh^{p+q}(R,k)^\Gamma\simeq
\coh^{p+q}(R\# k\Gamma,k) .
}
\end{equation}
Moreover, if $M$ is a finitely generated $R\# k\Gamma$-module,
there is a spectral sequence converging to the cohomology of $R$ with coefficients  in $M$:
\begin{equation}
\label{sseqmod}
\xymatrix
{\coh^{p+q}(\Gr_{(p)} R, M) \ar@{=>}[r] & \coh^{p+q}(R,M) ,
}
\end{equation}
also compatible with the action of $\Gamma$.

We wish to apply Lemma \ref{fingen} to the spectral sequence (\ref{sseq2}) for the
filtered algebra $R$.
In order to do so,  we must find some permanent  cycles.

The Hochschild cohomology of $R$ in degree 2, with trivial coefficients,
was studied in \cite{MW}:
There is a linearly independent set of
2-cocycles $\xi_{\alpha}$ on $R$, indexed by the positive roots $\alpha$
(\cite[Theorem 6.1.3]{MW}).
We will use the notation $\xi_{\alpha}$ in place of the notation $f_{\alpha}$
used there.
As shown in \cite{MW}, these
2-cocycles may be expressed as functions at the chain level in the following way.
Let ${\bf x}^{\bf a}$, ${\bf x}^{\bf b}$ denote arbitrary PBW basis elements
(\ref{PBW}) of $u({\mathcal D},0,0)$.
Let $\widetilde{\bf x}^{\bf a}$, $\widetilde{\bf x}^{\bf b}$
denote corresponding elements in the
infinite dimensional algebra $U({\mathcal D},0)$ arising from the section
of the quotient map $U({\mathcal D},0)\rightarrow u({\mathcal D},0,0)$
for which PBW basis elements are sent to PBW basis elements.
Then
\begin{equation}\label{xialpha}
 \xi_{\alpha}({\bf x}^{\bf a}\ot {\bf x}^{\bf b})=c_{\alpha}
\end{equation}
where $c_{\alpha}$ is the coefficient of $\widetilde{x}_{\alpha}^{N_{\alpha}}$
in the product $\widetilde{{\bf x}}^{{\bf a}}\cdot \widetilde{{\bf x}}^{{\bf b}}$, and
${\bf x}^{\bf a}$, ${\bf x}^{\bf b}$ range over all pairs of PBW basis elements.
A direct proof that these are 2-cocycles is in \cite[\S6.1]{MW};
alternatively the proof in Lemma \ref{permanent2} below, that analogous
functions $f_{\alpha}$ in higher degrees
are cocycles for $u({\mathcal D},\lambda,0)$, applies
with minor modifications to the $\xi_{\alpha}$ in this context.

We wish to relate these functions $\xi_{\alpha}$ to
elements  on the $E_1$-page of the spectral sequence (\ref{sseq1}), found  in the previous section.
Recall that $\beta_1, \ldots, \beta_r$ enumerate  the positive roots, and
the integer $p_i$ is defined in (\ref{pi}).
By (\ref{xialpha}) and Lemma \ref{filtration},
we have
$
  \xi_{\beta_i}\downarrow_{F_{p_i-1}(R\otimes R)} = 0  $ but
$\xi_{\beta_i}\downarrow_{F_{p_i}
  (R\otimes R)} \not = 0.
$
Note that the filtration on $R$ induces a filtration on $R^+$, and thus a
(decreasing) filtration on the complex $C^\bu$ defined in (\ref{cdot}),
given by $F^{p}C^n = \{f: (R^+)^{\otimes n} \to k \, | \,
f\downarrow_{F_{p-1}((R^+)^{\otimes n})} = 0 \}$.
We conclude that $\xi_{\beta_i} \in F^{p_i}C^2$ but $\xi_{\beta_i} \not \in F^{p_i+1}C^2$.   Denoting the corresponding cocycle by the same letter, we further conclude that
$\xi_{\beta_i} \in \im \{ \coh^2(F^{p_i}C^\bu) \to \coh^2(C^\bu)\} = F^{p_i}\coh^2(R,k)$,
but $ \xi_{\beta_i} \not \in \im \{ \coh^2(F^{p_{i}+1}C^\bu)
\to \coh^2(C^\bu)  \} =F^{p_i+1}\coh^2(R,k)$.
Hence, $\xi_{\beta_i}$ can be identified with the corresponding nontrivial homogeneous element
in the associated graded complex:
$$\widetilde \xi_{\beta_i} \in F^{p_i}\coh^2(R,k)/F^{p_i+1}\coh^2(R,k)
\simeq E_\infty^{p_i, 2-p_i}.$$
Since $\xi_{\beta_i} \in F^{p_i}C^2$ but $\xi_{\beta_i} \not \in F^{p_i+1}C^2$, it induces an element
$\bar \xi_{\beta_i} \in E_0^{p_i, 2-p_i} = F^{p_i}C^2/F^{p_i+1}C^2$.  Since $\bar \xi_{\beta_i}$  is induced by an actual cocycle in $C^\bu$,
it will be in the kernels of all the  differentials  of the spectral sequence.  Hence, the residue of $\bar \xi_{\beta_i}$ will be in the  $E_\infty$--term where it will be identified with the non-zero  element $\widetilde \xi_{\beta_i}$  since these classes are induced by the same  cocycle in $C^\bu$.
We conclude that $\bar \xi_{\beta_i} \in E^{p_i, 2-p_i}_0$, and, correspondingly, its image in $E^{p_i, 2-p_i}_1 \simeq \coh^2(\Gr R,k)$ which we denote by the same symbol,
is a permanent cycle.

 Note that the results of \cite{MW} apply
equally well to $S=\Gr R$ to yield similar cocycles $\hat \xi_\alpha$ via the formula
(\ref{xialpha}), for each positive root $\alpha$ in type $A_1\times\cdots\times A_1$.
Now $S=\Gr R$ has generators corresponding to the root vectors of $R$,
and we similarly identify elements in cohomology.
Thus each $\hat{\xi}_{\alpha}$, $\alpha$ a positive root in type
$A_1\times\cdots\times A_1$, may be relabeled $\hat{\xi}_{\beta_i}$ for some $i$
($1\leq i\leq r$), the indexing corresponding to $\beta_1,\ldots,\beta_r\in\Phi^+$.
Comparing the values of $\overline{\xi}_{\beta_i}$ and $\hat \xi_{\beta_i}$
on basis elements ${\bf x}^{\bf a}\ot {\bf x}^{\bf b}$ of $\Gr R\otimes \Gr R$ we  conclude
that they are the same function.
Hence  $\hat \xi_{\beta_i} \in   E_1^{p_i, 2-p_i}$ are permanent  cycles.

We wish to identify these elements $\hat \xi_{\beta_i} \in \coh^2(\Gr R,k)$
with  the cohomology classes $\xi_i$ in $\coh^*(S,k)$ constructed in
Section \ref{sec-step0}, as we may thus exploit the algebra structure
of $\coh^*(S,k)$ given in Theorem \ref{step0}.
We use the same symbols $x_{\beta_i}$ to denote the corresponding root
vectors in $R$ and in $S=\Gr R$, as this should cause no confusion.

\begin{lemma}
\label{perm}
For each $i$ ($1\leq i\leq r$),
the cohomology classes $\xi_i$ and $\hat \xi_{\beta_i}$ coincide as elements of
$\coh^2(\Gr R,k)$.
\end{lemma}

\begin{proof}
Let $K_{\bullet}$ be the chain complex defined in Section \ref{sec-step0},
a projective resolution of the trivial $\Gr R$-module $k$.
Elements $\xi_i\in\coh^2(\Gr R,k)$ and $\eta_i\in\coh^1(\Gr R,k)$ were
defined via the complex $K_{\bullet}$. We wish to identify
$\xi_i$ with elements of the chain complex $C^{\bullet}$ defined
in (\ref{cdot}), where $A=R$.
To this end we define maps $F_1, F_2$ making the following diagram
commute, where $S=\Gr R$:
$$
\begin{array}{ccccccccccc}
\cdots \! &\rightarrow & K_2 & \stackrel{d}{\longrightarrow}
  & K_1 & \stackrel{d}{\longrightarrow} & K_0 &
  \stackrel{\varepsilon}{\longrightarrow} & k & \rightarrow & \! 0\\
 & & \hspace{.1in}\downarrow F_2 & & \hspace{.1in}\downarrow F_1 & &
   \parallel &&\parallel &&\\
 \cdots \! &\! \rightarrow \! & S\ot (S^+)^{\ot 2} \! & \!
  \stackrel{\partial_2}{\longrightarrow}\!
  & \! S\ot S^+ \! &\! \stackrel{\partial_1}{\longrightarrow}
  \! & \! S &\stackrel
  {\varepsilon}{\longrightarrow}\! & \! k \! &\! \rightarrow & \! 0
\end{array}
$$
In this diagram, the maps $d$ are given in (\ref{dsum}), and
$\partial_i$ in (\ref{free-res}).
Let $\Phi(\cdots 1_i \cdots)$ denote the basis element of $K_1$ having
a 1 in the $i$th position, and 0 in all other positions, $\Phi(\cdots 1_i\cdots 1_j\cdots)$
(respectively $\Phi(\cdots 2_i\cdots)$) the basis element of $K_2$ having a 1 in the $i$th
and $j$th positions ($i\neq j$), and 0 in all other positions (respectively a 2 in the $i$th
position and 0 in all other positions).
Let
\begin{eqnarray*}
  F_1(\Phi(\cdots 1_i\cdots)) &=& 1\ot x_{\beta_i}, \\
  F_2(\Phi(\cdots 2_i\cdots)) &=& \sum_{a_i=0}^{N_i-2}
              x_{\beta_i}^{a_i}\ot x_{\beta_i}\ot x_{\beta_i}^{N_i-a_i-1}, \\
  F_2(\Phi(\cdots 1_i\cdots 1_j\cdots)) &=& 1\ot x_{\beta_i}\ot x_{\beta_j} - q_{ij}
             \ot x_{\beta_j}\ot x_{\beta_i} .
\end{eqnarray*}
Direct computations show that the
two nontrivial squares in the diagram above commute.
So $F_1,F_2$ extend to maps
$F_i : K_i\rightarrow
S\ot (S^+)^{\ot i}$, $i\geq 1$, providing a chain map $F_{\bu} : K_{\bu}\rightarrow S
\ot (S^+)^{\ot \bu }$,
thus inducing isomorphisms on cohomology.

We now verify that the maps $F_1, F_2$ make the desired identifications.
By $\xi_{\beta_i}$ we mean the function on the reduced bar complex,
 $
 \xi_{\beta_i}(1\ot {\bf x}^{\bf a}\ot {\bf x}^{\bf b})
   := \xi_{\beta_i}( {\bf x}^{\bf a}\ot {\bf x}^{\bf b})$,
defined in (\ref{xialpha}).
Then
\begin{eqnarray*}
  F_2^*(\xi_{\beta_i})(\Phi(\cdots 2_i \cdots)) & = & \xi_{\beta_i}
   (F_2(\Phi(\cdots 2_i\cdots)))\\
  &=& \xi_{\beta_i}(\sum_{a_i=0}^{N_i-2} x_{\beta_i}^{a_i}\ot x_{\beta_i}\ot x_{\beta_i}^{N_i-a_i-1})\\
  &=& \sum_{a_i=0}^{N_i-2} \varepsilon(x_{\beta_i}^{a_i})\xi_{\beta_i}
        (1\ot x_{\beta_i}\ot x_{\beta_i}^{N_i-a_i-1})\\
  &=& \xi_{\beta_i}(x_{\beta_i}\ot x_{\beta_i}^{N_i-1}) \ \ = \ \ 1.
\end{eqnarray*}
Further, it may be checked similarly that $F_2^*(\xi_{\beta_i})
(\Phi(\cdots 1_i\cdots 1_j\cdots ))=0$ for all $i,j$ and
$F_2^*(\xi_{\beta_i})(\Phi(\cdots 2_j\cdots))=0$ for all $j\neq i$.
Therefore $F_2^*(\xi_{\beta_i})$ is the dual function to
$\Phi(\cdots 2_i\cdots)$, which is precisely $\xi_i$.
\end{proof}

Similarly, we identify the elements $\eta_i$ defined in Section \ref{sec-step0}
with functions at the chain level in cohomology: Define
\begin{equation}\label{etabeta}
  \eta_{\alpha}({\bf x}^{\bf a})=\left\{
   \begin{array}{ll} 1, & \mbox{ if } {\bf x}^{\bf a}=x_{\alpha}\\
                        0, & \mbox{ otherwise}\end{array}\right..
\end{equation}
The functions $\eta_{\alpha}$ represent a basis of $\coh^1(S,k)\simeq \Hom_k(S^+/(S^+)^2,k)$.
Similarly functions corresponding to the {\em simple} roots $\beta_i$ only
represent a basis of $\coh^1(R,k)\simeq \Hom_k (R^+/(R^+)^2,k)$.
A computation shows that $F_1^*(\eta_{\beta_i})(\Phi(\cdots 1_j\cdots)) =\delta_{ij}$,
 so that $F_1^*(\eta_{\beta_i})$ is the dual function to
$\Phi(\cdots 1_i\cdots )$.
Therefore $\eta_i$ and $\eta_{\beta_i}$ coincide as elements of $\coh^1(S,k)$.

For each $\alpha\in\Phi^+$, let $M_{\alpha}$ be any positive integer for which
$\chi_{\alpha}^{M_{\alpha}}=\varepsilon$. (For example, let $M_{\alpha}$ be
the order of $\chi_{\alpha}$.)
Note that $\xi_{\alpha}^{M_{\alpha}}$ is $\Gamma$-invariant:
By (\ref{action}),
$g\cdot \xi_{\alpha}^{M_{\alpha}} =\chi_{\alpha}^{-M_{\alpha}N_{\alpha}}(g)
\xi_{\alpha}^{M_{\alpha}} =\xi_{\alpha}^{M_{\alpha}}$.

Recall the notation $u({\mathcal D},0,0)\simeq R\# k\Gamma$,
where $R={\mathcal B}(V)$ is the Nichols algebra.

\begin{lemma}\label{malpha}
The cohomology algebra $\coh^*(u({\mathcal D},0,0),k)$
is finitely generated over the subalgebra
generated by all $\xi_{\alpha}^{M_{\alpha}}$ ($\alpha\in\Phi^+$).
\end{lemma}

\begin{proof}
Let $\xymatrix{ E_1^{*,*} \ar@{=>}[r] & \coh^*(R,k)}$ be the
spectral sequence (\ref{sseq1}), and let
$B^{*,*}$ be the bigraded subalgebra  of $E_1^{*,*}$  generated by the elements $\xi_i$.
By Lemma  \ref{perm} and the discussion prior to it, $B^{*,*}$ consists of
permanent cycles.
Since  $\xi_i$ is  $\bar \xi_{\beta_i}$
by Lemma \ref{perm}, it is in bidegree
$(p_i, 2-p_i)$.
Let $A^{*,*}$ be the subalgebra of $B^{*,*}$ generated by $\overline{\xi}_{\alpha}^{M_{\alpha}}$
($\alpha\in\Phi^+$).
Then $A^{*,*}$ also consists of permanent cycles.
 Observe that
$A^{*,*}$ is a subalgebra of $\coh^*((\Gr R)\# k\Gamma,k)$, which is  graded commutative since
$(\Gr R)\# k\Gamma$ is a Hopf algebra.
Hence,  $A^{*,*}$ is commutative as    it is concentrated in even (total) degrees.
Finally, $A^{*,*}$ is Noetherian since it is a polynomial algebra in
the $\xi_{\alpha}^{M_{\alpha}}$.
We conclude that the  bigraded commutative algebra $A^{*,*}$ satisfies
the hypotheses of  Lemma~\ref{fingen}.  By Theorem~\ref{step0}, the algebra
$E_1^{*,*} \simeq \coh^*(\Gr R,k)$  is generated by $\xi_i$ and $\eta_i$ where
the generators  $\eta_i$ are nilpotent.
Hence, $E_1^{*,*}$ is a finitely generated module over $B^{*,*}$.
The latter is clearly a finitely generated  module over
$A^{*,*}$.
Hence  $E_1^{*,*}$ is a finitely generated module over $A^{*,*}$.
Lemma~\ref{fingen} implies that $\coh^*(R,k)$ is a Noetherian $\Tot(A^{*,*})$-module;
moreover, the action of $\Gamma$ on
$\coh^*(R,k)$ is compatible with the action on $A^{*,*}$ since the spectral sequence
is compatible with the action of $\Gamma$.
Therefore, $\coh^*(R\# k\Gamma,k) \simeq \coh^*(R,k)^\Gamma$ is a Noetherian
$\Tot(A^{*,*})$-module.
Since $\Tot(A^{*,*})$   is finitely generated, we conclude
that $\coh^*(R\# k\Gamma,k)$ is finitely generated.
\end{proof}

We immediately have the following theorem.
The second statement of the theorem
follows by a simple application of the second statement  of Lemma~\ref{fingen}.

\begin{thm}\label{step1}
The algebra $\coh^*(u({\mathcal D},0,0),k)$ is finitely generated.
If $M$ is a finitely generated  $u({\mathcal D},0,0)$-module,
then $\coh^*(u({\mathcal D},0,0),M)$ is a
finitely generated module over $\coh^*(u({\mathcal D},0,0),k)$.
\end{thm}

Thanks to Corollary \ref{cor:braided-graded}, we have some information
about the algebra structure of the cohomology ring of the Nichols algebra
$R={\mathcal B}(V)$:
Recall that $q_{\beta\alpha} = \chi_{\alpha}(g_{\beta})$ (see
(\ref{gbeta-chibeta})).
Compare the following result with the graded case, Theorem \ref{step0}.

\begin{thm}\label{relns}
The following relations hold in $\coh^*(R,k)$ for all
$\xi_{\alpha},\xi_{\beta}$ ($\alpha,\beta\in\Phi^+$) and
$\eta_{\alpha},\eta_{\beta}$ ($\alpha,\beta\in\Pi$):
$$
  \xi_{\alpha}\xi_{\beta} = q_{\beta\alpha}^{N_{\alpha}N_{\beta}}
  \xi_{\beta}\xi_{\alpha}, \ \ \
 \eta_{\alpha}\xi_{\beta} = q_{\beta\alpha}^{N_{\beta}}\xi_{\beta}\eta_{\alpha}, \mbox{ and }
   \ \ \ \eta_{\alpha}\eta_{\beta}=-q_{\beta\alpha}\eta_{\beta}\eta_{\alpha}.
$$
\end{thm}

\begin{proof}
As we shall see, this is a consequence of Corollary
\ref{cor:braided-graded}, since $R$ is a braided Hopf algebra in
$\YD$.

Note first that as a function on the Yetter-Drinfeld module
$R\in\YD$ the cocycle $\xi_{\alpha}$ is $\Gamma$-homogeneous of
degree $g_\alpha^{-N_\alpha}$ and spans a one-dimensional
$\Gamma$-module with character $\chi_\alpha^{-N_\alpha}$.
To see this in full detail, let us rephrase the definition
\eqref{xialpha} of $\xi_\alpha$ as follows.
One can write $U(\mathcal D,0)$ as a
Radford biproduct $U(\mathcal D,0)=\hat R\# k\Gamma$ with a
braided Hopf algebra $\hat R\in\YD$ such that $R$ is a quotient of
$\hat R$. Let $s\colon R\to\hat R$ be the section of the
surjection $\hat R\to R$ that maps PBW basis elements to PBW basis
elements. Note that $s$ is a map in $\YD$. Finally, let
$p_\alpha\colon\hat R\to k$ be the function projecting $r\in\hat
R$ to the coefficient in $r$ of the PBW basis element $\tilde
x_\alpha^{N_\alpha}$. Then $\xi_\alpha$ is by definition the
pullback of $p_\alpha$ under $R\ot R\xrightarrow{s\ot s}\hat
R\ot\hat R\xrightarrow{m}\hat R$. As a consequence, since
$p_\alpha$ is clearly $\Gamma$-homogeneous of degree
$g_\alpha^{-N_\alpha}$ and satisfies $g\cdot
p_\alpha=\chi_\alpha^{-N_\alpha}(g)p_\alpha$ (cf.\ (\ref{action})),
the statements on $\xi_\alpha$ follow.

On the other hand, it is easy to see from definition
(\ref{etabeta}) that $\eta_\alpha$ has degree $g_\alpha\inv$ and
spans a one-dimensional $\Gamma$-module with character
$\chi_\alpha\inv$.

Let us denote the opposite of multiplication in $\coh^*(R,k)$ by
$\zeta\circ\theta:=\theta\zeta$. According to
Corollary \ref{cor:braided-graded}, the opposite multiplication is braided
graded commutative. In particular
\begin{align*}
  \xi_\alpha\xi_\beta=\xi_\beta\circ\xi_\alpha
    &=\left(g_\beta^{-N_\beta}\cdot\xi_\alpha\right)\circ\xi_\beta\\
    &=\chi_\alpha^{-N_\alpha}\left(g_\beta^{-N_\beta}\right)\xi_\alpha\circ\xi_\beta\\
    &=\chi_\alpha\left(g_\beta\right)^{N_\beta N_\alpha}\xi_\alpha\circ\xi_\beta\\
    &=q_{\beta\alpha}^{N_\beta N_\alpha}\xi_\alpha\circ\xi_\beta
    =q_{\beta\alpha}^{N_\beta N_\alpha}\xi_\beta\xi_\alpha
\end{align*}
and \[  \eta_\alpha\xi_\beta=\xi_\beta\eta_\alpha
    =\left(g_\beta^{-N_\beta}\cdot\eta_\alpha\right)\circ\xi_\beta
    =\chi_{\alpha}\inv\left(g_\beta^{-N_\beta}\right)\eta_\alpha\circ\xi_\beta
    =q_{\beta\alpha}^{N_\beta}\xi_\beta\eta_\alpha\]
as well as finally
\[\eta_\alpha\eta_\beta=\eta_\beta\circ\eta_\alpha
    =-\left(g_\beta\inv\cdot\eta_\alpha\right)\circ\eta_\beta
    =-\chi_\alpha\inv\left(g_\beta\inv\right)\eta_\alpha\circ\eta_\beta
    =-q_{\beta\alpha}\eta_\beta\eta_\alpha.\]
\end{proof}

We give a corollary in a special case.
It generalizes  \cite[Theorem 2.5(i)]{GK}.
Recall that $\Pi$ denotes a set of simple roots in the root system $\Phi$.

\begin{cor}\label{no-odd}
Assume there are no $\Gamma$-invariants in $\coh^*(\Gr R,k)$ of the form
$\xi_{\beta_1}^{b_1}\cdots\xi_{\beta_r}^{b_r}\eta_{\beta_1}^{c_1}\cdots \eta_{\beta_r}
^{c_r}$ for which $c_1+\cdots +c_r$ is odd. Then
 $$\coh^*(u({\mathcal D},0,0), k)\simeq k\langle\xi_{\alpha},\eta_{\beta}\mid
   \alpha\in\Phi^+, \ \beta\in \Pi\rangle^{\Gamma},$$
with the relations of Theorem \ref{relns}, $\deg(\eta_{\beta})=1$,
$\deg(\xi_{\alpha})=2$.
If there are no $\Gamma$-invariants with $c_1+\cdots +c_r\neq 0$,
then $\coh^*(R\# k\Gamma,k)\simeq k[\xi_{\alpha}\mid \alpha\in\Phi^+]^{\Gamma}$.
\end{cor}

\begin{proof}
The hypothesis of the first statement implies that $\coh^i(\Gr R,k)^{\Gamma}=0$
for all odd integers $i$.
Thus on the $E_1$-page, every other diagonal is 0.
This implies $E_1 = E_{\infty}$.
It follows that, as a vector space, $\coh^*(u({\mathcal D},0,0),k)$ is exactly as stated.
The algebra structure is a consequence of Theorem \ref{relns}
and the fact that the cohomology of a Hopf algebra is graded
commutative.
The hypothesis of the last statement implies further that the $\Gamma$-invariant
subalgebra of $\coh^*(\Gr R,k)$ is spanned by elements of the form
$\xi_{\beta_1}^{b_1}\cdots \xi_{\beta_r}^{b_r}$.
By graded commutativity of the cohomology ring of a Hopf algebra
and the relations of Theorem \ref{relns}, $\coh^*(R\# k\Gamma, k)$ may
be identified with the $\Gamma$-invariant subalgebra of a polynomial ring
in variables $\xi_{\alpha}$, with corresponding $\Gamma$-action.
\end{proof}

\begin{rem}
{\em Assume the hypotheses of Corollary \ref{no-odd}, and that
$q_{\alpha\alpha}\neq -1$ for all $\alpha\in\Pi$.
Then $\eta_{\alpha}^2=0$ for all $\alpha\in\Pi$, and it follows that
the maximal ideal spectrum of $\coh^*(u({\mathcal D},0,0),k)$ is
$\Spec k[\xi_{\alpha}\mid\alpha\in\Phi^+]^{\Gamma} \simeq
\Spec k[\xi_{\alpha}\mid\alpha\in\Phi^+]/\Gamma$.}
\end{rem}

We give an example to show that cohomology may in fact be nonzero
in odd degrees,
in contrast to that of the small quantum groups of \cite{GK}.
This complicates any determination of the explicit structure of
cohomology in general.  The  simplest example occurs  in type
$A_1\times A_1 \times A_1$, where there can exist
a nonzero cycle  in degree $3$.

\begin{example}
{\em
Let $\Gamma = \Z/ \ell \Z$   with generator $g$.
Let $q$ be a primitive $\ell$th root of unity.
Let $g_1 = g_2 =g_3= g$, and  choose $\chi_1,\chi_2,\chi_3$ so that  the matrix $(q_{ij}) = (\chi_j(g_i))$ is
$$
  \left(\begin{array}{ccc}
   q&q^{-1}&1\\
   q&q^{-2}&q\\
1&q^{-1}&q\\
    \end{array}\right).
$$
Let $u({\mathcal D},0,0)=R \# k\Gamma$
be the  pointed  Hopf  algebra  of type $A_1\times A_1 \times A_1$ defined by this data.
Let $\eta_1, \eta_2, \eta_3$ represent elements of $\coh^1(R,k)$ as defined
by (\ref{etabeta}).
The action of $\Gamma$ is described in (\ref{action}): $ g_i\cdot \eta_j =  q_{ij}^{-1} \eta_j$.
Since the  product of all entries in any given row of the
matrix $(q_{ij})$ is 1, we conclude that $\eta_1\eta_2\eta_3$ is invariant  under $\Gamma$.
Hence,   it is a nontrivial cocycle in $\coh^3(u({\mathcal D},0,0), k)$. }
\end{example}

We also give an example in type $A_2\times A_1$ to illustrate, in particular,
that the methods  employed in  \cite{GK}  do not transfer to our more general setting.
That is, for an arbitrary (coradically graded) pointed  Hopf  algebra,  the first  page
of the  spectral sequence (\ref{sseq2}) can have nontrivial elements in odd degrees.
In the special case of a small quantum group
(with some restrictions on the order  of the root of unity), it is shown in \cite{GK}
that this does not happen.

\begin{example}\label{odd}
{\em
Let $\Gamma = \Z/\ell \Z \times \Z/\ell \Z$, with generators
$g_1,g_2$. Let $q$ be a primitive $\ell$th root of unity ($\ell$ odd), and let $\mathcal D$  be of type $A_2\times A_1$ so that the Cartan matrix is
$$
  \left(\begin{array}{rrr}
   2&-1&0\\
   -1&2&0\\
   0&0&2\end{array}\right).
$$
Let
$\chi_1,\chi_2$ be as for $u_q(sl_3)^+$, that is
$$
\begin{array}{rclrcl}
  \chi_1(g_1) & = & q^2,  \ \ \ & \ \ \ \chi_1(g_2) &=& q^{-1},\\
  \chi_2(g_1) &=& q^{-1}, & \chi_2(g_2)&=&q^2.
\end{array}
$$
Now let $g_3:= g_1g_2$ and $\chi_3:= \chi_1^{-1}\chi_2^{-1}$.
Then $\chi_3(g_3)=q^{-2}\neq 1$ and the Cartan condition holds,
for example
  $\chi_3(g_1)\chi_1(g_3)= q^{-1} q = 1 = \chi_1(g_1)^{a_{13}}$
since $a_{13} =0$.
Let $R={\mathcal B}(V)$, the Nichols algebra defined from this
data.

The root vector corresponding to the nonsimple positive root is
$x_{12} = [x_1,x_2]_c=x_1x_2-q^{-1}x_2x_1$.
The relations among the root vectors other than $x_3$ are now
$$
  x_2x_1 = qx_1x_2 - qx_{12}, \ \ \ x_{12}x_1=q^{-1}x_1x_{12}, \ \ \
   x_2x_{12}=q^{-1}x_{12}x_2.
$$
The associated graded algebra, which is of type $A_1\times A_1\times
A_1\times A_1$,
thus has relations (excluding those
involving $x_3$ that do not change):
$$
  x_2x_1 = qx_1x_2, \ \ \ x_{12}x_1=q^{-1}x_1x_{12}, \ \ \
   x_2x_{12}=q^{-1}x_{12}x_2.
$$
Thus the cohomology of the associated graded algebra $\Gr R \# k\Gamma$
is the subalgebra of $\Gamma$-invariants of an algebra with generators
$$
  \xi_1, \ \xi_{12}, \ \xi_2, \ \xi_3, \ \eta_1, \ \eta_{12}, \ \eta_2, \ \eta_3
$$
(see Theorem \ref{step0}).
The element $\eta_1\eta_2\eta_3$ is $\Gamma$-invariant
since $g\cdot\eta_i = \chi_i(g)^{-1} \eta_i$ for all $i$,
and $\chi_3=\chi_1^{-1}\chi_2^{-1}$.
Therefore in the spectral sequence relating the cohomology
of $\Gr R\# k\Gamma$ to that of $R\# k\Gamma$, there are some
odd degree elements on the $E_1$ page.
}\end{example}

\section{Finite dimensional pointed Hopf algebras} \label{sec-step2}

Let $u({\mathcal D},\lambda,\mu)$ be one of the finite dimensional
pointed Hopf algebras
from the Andruskiewitsch-Schneider classification \cite{AS06}, as described
in Section \ref{pointed}.
With respect to the coradical filtration, its associated graded algebra is
$\Gr u({\mathcal D},\lambda,\mu)\simeq u({\mathcal D},0,0)\simeq R\# k\Gamma$ where
$R={\mathcal B}(V)$ is a Nichols algebra, also described in Section \ref{pointed}.
By Theorem \ref{step1},
$\coh^*(u({\mathcal D},0,0), k)$
is finitely generated.
Below we
 apply the  spectral sequence for a filtered algebra, employing this new choice of filtration:

\begin{equation}\label{sseq10}
\xymatrix
{E^{p,q}_1 = \coh^{p+q}(u({\mathcal D},0,0), k) \ar@{=>}[r] & \coh^{p+q}(
u({\mathcal {D}},\lambda,\mu),k) .
}
\end{equation}

As a consequence of
the following lemma, we may assume without loss of generality that
all root vector relations (\ref{root-relns}) are trivial.

\begin{lemma}\label{mu-zero}
(i) For all ${\mathcal D}$, $\lambda$, $\mu$, there is an isomorphism of graded algebras,
$$\coh^*(u({\mathcal D}, \lambda,\mu),k)\simeq \coh^*(u({\mathcal D}, \lambda,0),k).$$

(ii) Let $M$ be a finitely generated $u({\mathcal D}, \lambda,\mu)$-module.
There exists a finitely generated $u({\mathcal D}, \lambda,0)$-module
$\widetilde M$ and an isomorphism of $\coh^*(u({\mathcal D}, \lambda,\mu),k)$-modules
$$
\coh^*(u({\mathcal D}, \lambda,\mu), M)\simeq \coh^*(u({\mathcal D}, \lambda,0),\widetilde M),
$$
where the action of $\coh^*(u({\mathcal D}, \lambda,\mu),k)$ on $\coh^*(u({\mathcal D}, \lambda,0),\widetilde M)$  is via the isomorphism
of graded algebras in (i).

\end{lemma}

\begin{proof} Define the subset $I \subset \Gamma$ as follows: $g \in I$  if and  only if there exists a root
vector relation $x_\alpha^{N_\alpha} - u_\alpha(\mu)$ (see  (\ref{root-relns})), such that
the element $u_\alpha(\mu)$ of $k\Gamma$
has a nonzero coefficient of $g$ when written as a linear combination of group elements.
If $x_\alpha^{N_\alpha} - u_\alpha(\mu)$ is a nontrivial root vector relation, then necessarily
$x_{\alpha}^{N_{\alpha}}$ commutes with elements of $\Gamma$, implying that $\chi_{\alpha}^{N_{\alpha}}=\varepsilon$.
{F}rom this and (\ref{Ncentral}), we see that
$x_{\alpha}^{N_{\alpha}}$ is central in $u({\mathcal D},\lambda,\mu)$.
Since $x_{\alpha}^{N_{\alpha}} = u_\alpha(\mu)$, this element  of the group ring $k\Gamma$ must
also be central. Since $\Gamma$ acts diagonally,
each group element involved in $u_\alpha(\mu)$ is necessarily central
in $u({\mathcal D},\lambda,\mu)$ as well.

Let $Z=k\langle I\rangle$, a subalgebra
of $u({\mathcal D}, \lambda,\mu)$, and let $\overline{u}=u({\mathcal D}, \lambda,\mu)/(g -1\mid g\in I)$.
We have a sequence of algebras (see \cite[\S5.2]{GK}):
$$
   Z\rightarrow u({\mathcal D},\lambda,\mu) \rightarrow
      u({\mathcal D},\lambda,\mu)/\! /Z\simeq\overline{u}.
$$
Hence, there is a multiplicative spectral sequence
$$\coh^p(\overline{u}, \coh^q(Z, k))  \implies \coh^{p+q}(u({\mathcal D},\lambda,\mu),k).
$$
Since
the characteristic of $k$ does not divide the order of the group, we have $\coh^{q>0}(Z, k) = 0$. Thus
the spectral sequence collapses, and we get an isomorphism  of graded algebras
\begin{equation}
\label{spec_seq}
\coh^*(\overline{u}, k) \simeq \coh^*(u({\mathcal D},\lambda,\mu), k).
\end{equation}
Similarly, if $M$ is any  $u({\mathcal D},\lambda,\mu)$-module, then we have a spectral sequence of $\coh^*(\overline{u},k)$-modules
$$\coh^p(\overline{u}, \coh^q(Z, M))  \implies \coh^{p+q}(u({\mathcal D},\lambda,\mu),M).
$$
The spectral sequence collapses, and we get an isomorphism
\begin{equation}
\label{spec_seq_mod}
\coh^*(\overline{u}, M^Z) \simeq \coh^*(u({\mathcal D},\lambda,\mu), M)
\end{equation}
which respects  the action of $\coh^*(u({\mathcal D},\lambda,\mu),k)$, where $\coh^*(u({\mathcal D},\lambda,\mu),k)$ acts
on the left side via the isomorphism (\ref{spec_seq}).

Note that $Z =k\langle I\rangle$ is also a central subalgebra  of $u({\mathcal D},\lambda,0)$. Arguing exactly as above,  we get an isomorphism of graded algebras
$$
\coh^*(\overline{u}, k) \simeq \coh^*(u({\mathcal D},\lambda,0), k),
$$
which implies (i).

Let $M$ be a $u({\mathcal D},\lambda, \mu)$-module. Let $\widetilde M$ be a $u({\mathcal D},\lambda,0)$-module  which we get by inflating
the $\overline{u}$-module  $M^Z$ via the projection
$\xymatrix{u({\mathcal D},\lambda,0) \ar@{->>}[r] &\overline u }$.   Since $(\widetilde M)^Z \simeq M^Z$ by construction, we get an isomorphism  of
$u({\mathcal D},\lambda,0)$-modules
$$
\coh^*(\overline u, M^Z) \simeq \coh^*(u({\mathcal D},\lambda,0), \widetilde M)
$$
using another spectral sequence argument.
Combining with the isomorphism (\ref{spec_seq_mod}), we get (ii).

\end{proof}

By Lemma \ref{mu-zero},
it suffices to work with the cohomology of $u({\mathcal{D}},\lambda,0)$,
in which all the root vectors are nilpotent.
In this case we define some permanent cycles:
As before, for each $\alpha\in\Phi^+$, let $M_{\alpha}$ be  any positive integer for which
$\chi_{\alpha}^{M_{\alpha}}=\varepsilon$ (for example, take $M_{\alpha}$ to
be the order of $\chi_{\alpha}$).
If $\alpha=\alpha_i$ is simple, then $\chi_{\alpha}(g_{\alpha})$ has order $N_{\alpha}$, and so
$N_{\alpha}$ divides $M_{\alpha}$.

We  previously identified an element $\xi_{\alpha}$  of $\coh^2(R,k)$,
where $R={\mathcal B}(V)$.
Now ${\mathcal B}(V)$ is no longer a subalgebra of $u({\mathcal D},\lambda,0)$
in general, due to the potential existence of
nontrivial linking relations, but we will show that still there
is an element analogous to $\xi_{\alpha}^{M_{\alpha}}$
in $\coh^{2M_{\alpha}}(u({\mathcal D}, \lambda, 0),k)$.
For simplicity, let $U=U({\mathcal D},\lambda)$
and $u=u({\mathcal D},\lambda,0)\simeq U({\mathcal D},\lambda)
/(x_{\alpha}^{N_{\alpha}}\mid \alpha\in \Phi^+)$ (see Section \ref{pointed}).
Let $U^+$, $u^+$ denote the augmentation ideals of $U$, $u$.

We will use a similar construction as in Section 4 of \cite{MW},
defining functions as elements of the bar complex (\ref{cdot}).
For each $\alpha\in\Phi^+$,
define a $k$-linear function $\widetilde{f}_{\alpha}: (U^+)^{2M_{\alpha}}\rightarrow k$
by  first
letting $r_1,\ldots,r_{2M_{\alpha}}$ be PBW basis elements (\ref{PBW}) and
requiring
$$
  \widetilde{f}_{\alpha}(r_1\ot\cdots\ot r_{2M_{\alpha}}) = \gamma_{12}\gamma_{34}\cdots
  \gamma_{2M_{\alpha}-1,2M_{\alpha}}
$$
where $\gamma_{ij}$ is the coefficient of $x_{\alpha}^{N_{\alpha}}$
in the product $r_ir_j$ as a linear combination of
PBW basis elements.
Now define $\widetilde{f}_{\alpha}$ to be 0 whenever a tensor factor is in
$k\Gamma\cap \ker\varepsilon$, and
$$
  \widetilde{f}_{\alpha}(r_1g_1\ot\cdots\ot r_{2M_{\alpha}}g_{2M_{\alpha}})=
  \widetilde{f}_{\alpha}(r_1\ot {}^{g_1}r_2\ot\cdots\ot {}^{g_1\cdots g_{2M_{\alpha}-1}}r_{2M_{\alpha}})
$$
for all $g_1,\ldots,g_{2M_{\alpha}}\in\Gamma$.
It follows from the definition of $\widetilde{f}_{\alpha}$ and the fact
that $\chi_{\alpha}^{M_{\alpha}}=\varepsilon$ that $\widetilde{f}_{\alpha}$ is
$\Gamma$-invariant.
We will show that $\widetilde{f}_{\alpha}$ factors through the quotient $u^+$ of $U^+$ to give
a map $f_{\alpha}: (u^+)^{\ot 2M_{\alpha}}\rightarrow k$.
Precisely, it suffices to show that $\widetilde{f}_{\alpha}(r_1\ot\cdots\ot r_{2M_{\alpha}})=0$
whenever one of $r_1,\ldots,r_{2M_{\alpha}}$ is in the kernel of the quotient
map $\pi : U({\mathcal D},\lambda)\rightarrow u({\mathcal D},\lambda,0)$.
Suppose $r_i\in\ker \pi$, that is
$r_i=x_{\beta_1}^{a_1}\cdots x_{\beta_r}^{a_r}$ and for some $j$,
$a_j\geq N_{\beta_j}$.
Since $x_{\beta_j}^{N_{\beta_j}}$ is braided-central, $r_i$ is a scalar multiple
of $x_{\beta_j}^{N_{\beta_j}} x_{\beta_1}^{b_1}\cdots x_{\beta_r}^{b_r}$ for some
$b_1,\ldots,b_r$.
Now $\widetilde{f}_{\alpha}(r_1\ot\cdots\ot r_{2M_{\alpha}})$ is the product of the
coefficients of $x_{\alpha}^{N_{\alpha}}$ in $r_1r_2, \ldots, r_{2M_{\alpha}-1}
r_{2M_{\alpha}}$.
However, the coefficient of $x_{\alpha}^{N_{\alpha}}$ in each of $r_{i-1}r_i$
and $r_ir_{i+1}$ is 0:
If $\alpha =\beta_i$, then since $r_{i-1}, r_{i+1}\in U^+$,
this product cannot have a nonzero coefficient for $x_{\alpha}^{N_{\alpha}}$.
If $\alpha\neq \beta_i$, the same is true since $x_{\beta_j}^{N_{\beta_j}}$
is a factor of $r_{i-1}r_i$ and of $r_ir_{i+1}$.
Therefore $\widetilde{f}_{\alpha}$ factors to give a linear map
$f_{\alpha}: (u^+)^{2M_{\alpha}}\rightarrow k$.
In calculations, we define $f_{\alpha}$ via $\widetilde{f}_{\alpha}$ and
a choice of section of the quotient map $\pi: U\rightarrow u$.

\begin{lemma}\label{permanent2}
For each $\alpha\in\Phi^+$, $f_{\alpha}$ is a cocycle.
The $f_{\alpha}$ ($\alpha\in\Phi^+$) represent a linearly independent
subset of $\coh^*(u({\mathcal D},\lambda,0),k)$.
\end{lemma}

\begin{proof}
We first verify that $\widetilde{f}_{\alpha}$ is a cocycle on $U$:
Let $r_0,\ldots,r_{2M_{\alpha}}\in U^+$, of positive degree.
Then
$$
   d(\widetilde{f}_{\alpha})(r_0\ot\cdots\ot r_{2M_{\alpha}})=
  \sum_{i=0}^{2M_{\alpha}-1} (-1)^{i+1} \widetilde{f}_{\alpha}(r_0\ot\cdots\ot r_ir_{i+1}\ot\cdots
  \ot r_{2M_{\alpha}}).
$$
By definition of $\widetilde{f}_{\alpha}$, note that the first two terms cancel:
$$
  \widetilde{f}_{\alpha}(r_0r_1\ot r_2\ot\cdots\ot r_{2M_{\alpha}}) = \widetilde{f}_{\alpha}
 (r_0\ot r_1r_2\ot\cdots\ot r_{2M_{\alpha}}),
$$
and similarly for all other terms, so $d(\widetilde{f}_{\alpha})(r_0\ot\cdots\ot r_{2M_{\alpha}})
=0$. A similar calculation shows that
$d(\widetilde{f}_{\alpha})(r_0g_0\ot\cdots\ot r_{2M_{\alpha}-1}g_{2M_{\alpha}-1})=0$
for all $g_0,\ldots,g_{2M_{\alpha}-1}\in \Gamma$.
If there is an element of $k\Gamma\cap \ker\varepsilon$ in one of
the factors,  we obtain 0 as well by the definition of $\widetilde{f}_{\alpha}$,
a similar calculation.

Now we verify that $f_{\alpha}$ is a cocycle on the quotient $u$ of $U$:
Let $r_0,\ldots,r_{2M_{\alpha}}\in u^+$. Again we have
$$
   d(f_{\alpha})(r_0\ot\cdots\ot r_{2M_{\alpha}})=
  \sum_{i=0}^{2M_{\alpha}-1} (-1)^{i+1} f_{\alpha}(r_0\ot\cdots\ot r_ir_{i+1}\ot\cdots
  \ot r_{2M_{\alpha}}).
$$
We will show that
$
  f_{\alpha}(r_0r_1\ot r_2\ot\cdots\ot r_{2M_{\alpha}}) =
  f_{\alpha}(r_0\ot r_1r_2\ot\cdots\ot r_{2M_{\alpha}}),
$
   and similarly for the other terms.
Let $\widetilde{r}_i$ denote the element of $U$ corresponding to $r_i$ under
a chosen section of the quotient map $\pi: U\rightarrow u$.
Note that
$\widetilde{r}_0\cdot \widetilde{r}_1 = \widetilde{r_0r_1} + y$ and
$\widetilde{r}_1\cdot \widetilde{r}_2 = \widetilde{r_1r_2} + z$
for some $y,z\in\ker\pi$.
So
\begin{eqnarray*}
  f_{\alpha}(r_0r_1\ot r_2\ot \cdots \ot r_{2M_{\alpha}}) &=& \widetilde{f}_{\alpha}
    (\widetilde{r_0r_1} \ot \widetilde{r}_2\ot \cdots \ot \widetilde{r}_{2M_{\alpha}})\\
  &=&  \widetilde{f}_{\alpha} ((\widetilde{r}_0\cdot \widetilde{r}_1-y)
         \ot \widetilde{r}_2\ot \cdots \ot \widetilde{r}_{2M_{\alpha}})\\
 &=&  \widetilde{f}_{\alpha} (\widetilde{r}_0\cdot \widetilde{r}_1
         \ot \widetilde{r}_2\ot \cdots \ot \widetilde{r}_{2M_{\alpha}})\\
 &=&  \widetilde{f}_{\alpha} (\widetilde{r}_0 \ot\widetilde{r}_1\cdot
          \widetilde{r}_2\ot \cdots \ot \widetilde{r}_{2M_{\alpha}})\\
 &=&  \widetilde{f}_{\alpha} (\widetilde{r}_0 \ot(\widetilde{r}_1\cdot
          \widetilde{r}_2 + z)\ot \cdots \ot \widetilde{r}_{2M_{\alpha}})\\
 &=&  \widetilde{f}_{\alpha} (\widetilde{r}_0 \ot\widetilde{r_1r_2}
          \ot \cdots \ot \widetilde{r}_{2M_{\alpha}})\\
 &=& f_{\alpha}(r_0\ot r_1r_2\ot\cdots\ot r_{2M_{\alpha}}).
\end{eqnarray*}
Other computations for this case are similar to those for $U$.
Thus $f_{\alpha}$ is a cocycle on $u$.

We prove that in a given degree, the $f_{\alpha}$
 in that degree represent a linearly independent set in cohomology:
Suppose $\sum_{\alpha}c_{\alpha}f_{\alpha}=\partial h$
for some scalars $c_{\alpha}$ and linear map $h$.
Then for each $\alpha$,
\begin{eqnarray*}
  c_{\alpha} &=& (\sum c_{\alpha}f_{\alpha})
   (x_{\alpha}\ot x_{\alpha}^{N_{\alpha}-1}\ot\cdots\ot x_{\alpha}
   \ot x_{\alpha}^{N_{\alpha}-1})\\
  &=& \partial h (x_{\alpha}\ot x_{\alpha}^{N_{\alpha}-1}\ot\cdots\ot x_{\alpha}
   \ot x_{\alpha}^{N_{\alpha}-1})\\
  &=& -h(x_{\alpha}^{N_{\alpha}}\ot x_{\alpha}\ot\cdots\ot x_{\alpha}
  ^{N_{\alpha}-1}) + \cdots - h(x_{\alpha}\ot\cdots
    \ot x_{\alpha}^{N_{\alpha}-1}\ot x_{\alpha}^{N_{\alpha}})\ \ = \ \ 0
\end{eqnarray*}
since $x_{\alpha}^{N_{\alpha}}=0$ in $u=u({\mathcal D},\lambda,0)$.
\end{proof}

These functions $f_{\alpha}$
correspond to their counterparts $\xi_{\alpha}^{M_{\alpha}}$ defined on
$u({\mathcal D},0,0)$,
in the $E_1$-page of the spectral sequence (\ref{sseq10}),
by observing what they do as functions at the level of chain complexes
(\ref{cdot}).

\begin{thm}\label{mainthm}
The algebra $\coh^*(u({\mathcal D},\lambda,\mu),k)$ is finitely generated.
If $M$ is a finitely generated $u({\mathcal D},\lambda,\mu)$-module,
then $\coh^*(u({\mathcal D},\lambda,\mu),M)$ is a finitely generated module over
$\coh^*(u({\mathcal D},\lambda,\mu),k)$.
\end{thm}

\begin{proof}
By Lemma \ref{mu-zero}, it suffices to prove the statements in the case
$\mu=0$. We have $E_1^{*,*}\simeq \coh^*(u({\mathcal D},0,0),k)$, where
$u({\mathcal D},0,0)\simeq R\# k\Gamma$, $R={\mathcal B}(V)$.
By Lemma \ref{malpha}, $E_1^{*,*}$
is finitely generated over its subalgebra that is generated by all
$\xi_{\alpha}^{M_{\alpha}}$.
By Lemma \ref{permanent2} and the above remarks,
each $\xi_{\alpha}^{M_{\alpha}}$ is a permanent
cycle, corresponding to the cocycle $f_{\alpha}$ on $u({\mathcal D},\lambda,\mu)$.
The rest of the proof is an application of
Lemma \ref{fingen}, where $A^{*,*}$ is the subalgebra of $E_1^{*,*}$ generated by
the $\xi_{\alpha}^{M_{\alpha}}$ ($\alpha\in\Phi^+$), similar to the
proof of Lemma \ref{malpha} and Theorem \ref{step1}.
\end{proof}

\begin{cor}
The Hochschild cohomology ring $\coh^*(u({\mathcal D},\lambda,\mu),
u({\mathcal D},\lambda,\mu))$ is finitely generated.
\end{cor}

\begin{proof}
Apply Theorem \ref{mainthm} to the finitely generated
$u({\mathcal D},\lambda,\mu)$-module $u({\mathcal D},\lambda,\mu)$,
under the adjoint action.
By \cite[Prop.\ 5.6]{GK}, this is isomorphic to the Hochschild
cohomology ring of $u({\mathcal D},\lambda,\mu)$.
\end{proof}

In the special case of a small quantum group, we obtain the following finite generation result (cf. \cite[Thm.\ 1.3.4]{BNPP}).

\begin{cor}  Let $u_q({\mathfrak{g}})$ be a quantized restricted enveloping algebra such that
the order $\ell$ of $q$ is odd and prime to 3 if $\mathfrak g$ is of type $G_2$.
Then $\coh^*(u_q({\mathfrak{g}}),k)$  is a finitely generated algebra.  Moreover, for any finitely generated
$u_q({\mathfrak{g}})$-module  $M$,  $\coh^*(u_q({\mathfrak{g}}),M)$ is a finitely generated module over $\coh^*(u_q({\mathfrak{g}}),k)$.
 \end{cor}

\begin{remark}{\em
The restrictions on $\ell$ correspond to the assumptions (\ref{assumptions}).
However our techniques and results are more general: The restrictions are used
in the classification of Andruskiewitsch and Schneider, but not in our arguments.
We need only the filtration lemma of De Concini and Kac \cite[Lemma 1.7]{DCK}
as generalized to our setting (Lemma \ref{filtration}) to guarantee existence of
the needed spectral sequences.
Our results should hold for all small quantum groups $u_q({\mathfrak{g}})$
having such a filtration, including those at even roots of unity $q$ for
which $q^{2d}\neq 1$ (see \cite{beck94} for the general theory at even roots
of unity).
}
\end{remark}

We illustrate the connection between our results and those in
\cite{GK} with a small example.
Our structure results are not as strong, however our finite generation
result is much more general.

\begin{example}{\em
As an algebra, $u_q(sl_2)$ is generated by
$E,F,K$,
with relations $K^{\ell}=1$, $E^{\ell}=0$, $F^{\ell}=0$,
$KEK^{-1}=q^2E$, $KFK^{-1}=q^{-1}F$, and
\begin{equation}\label{lk}
EF-FE = \frac{K-K^{-1}}{q-q^{-1}},
\end{equation}
where $q$ is a primitive $\ell$th root of 1, $\ell >2$.
Consider the coradical filtration on $u_q(sl_2)$, in which
$\deg(K)=0$, $\deg(E)=\deg(F)=1$.
Note that $\Gr u_q(sl_2)$ is generated
by $E,F,K$, with all relations being the same except that (\ref{lk})
is replaced by $EF-FE=0$. This is an algebra of the type
featured in Section~\ref{sec-step0}:
$\coh^*(\Gr u_q(sl_2),k)\simeq
k[\xi_1,\xi_2,\eta_1\eta_2]/((\eta_1\eta_2)^2)$, since these are the
invariants, under the action of $\Gamma = \langle K\rangle$,
of the cohomology of the subalgebra of $\Gr u_q(sl_2)$ generated by $E,F$.
By \cite{GK},
$$
  \coh^*(u_q(sl_2),k)\simeq k[\alpha,\beta,\gamma]/(\alpha\beta+\gamma^2),
$$
the coordinate algebra of the nilpotent cone of $sl_2$.
Identify $\alpha\sim \xi_1$,
$\beta\sim \xi_2$, $\gamma\sim \eta_1\eta_2$: Then as maps, $\deg(\alpha)=\ell$,
$\deg(\beta)=\ell$, $\deg(\gamma)=2$, so $\alpha\beta+\gamma^2=0$
will imply $\gamma^2=0$ in the associated graded algebra, as expected.
}\end{example}

\end{document}